\documentclass{article}%
\usepackage{amsmath}
\usepackage{amsfonts}
\usepackage{amssymb}
\usepackage{graphicx}%
\providecommand{\U}[1]{\protect\rule{.1in}{.1in}}
\newtheorem{theorem}{Theorem}

\newtheorem{corollary}[theorem]{Corollary}

\newtheorem{lemma}[theorem]{Lemma}

\newenvironment{proof}[1][Proof]{\noindent\textbf{#1.} }{\ \rule{0.5em}{0.5em}}
\begin{document}

\title{Attractors for weak and strong solutions of the three-dimensional
Navier-Stokes equations with damping}
\author{Daniel Pardo$^{1}$, Jos\'{e} Valero$^{2}$, $\ $\'{A}ngel Gim\'{e}nez$^{2}$\\{\small ${}$Universidad Miguel Hernandez de Elche, Centro de Investigaci\'{o}n
Operativa,}\\{\small 03202} {\small Elche (Alicante), Spain.}\\{\small E.mails: daniel.pardo@alu.umh.es;\ jvalero@umh.es; a.gimenez@umh.es}}
\date{}
\maketitle

\begin{abstract}
In this paper we obtain the existence of global attractors for the dynamical
systems generated by weak solution of the three-dimensional Navier-Stokes
equations with damping.

We consider two cases, depending on the values of the parameters $\beta
,\alpha$ controlling the damping term and the viscosity $\mu$. First, for
$\beta\geq3$ we define a multivalued dynamical systems and prove the existence
of the global attractor as well. Second, for either $\beta>3$ or
$\beta=3,\ 4\alpha\mu>1$ the weak solutions are unique and we prove that the
global attractor for the corresponding semigroup is more regular. Also, we
prove in this case that it is the global attractor for the semigroup generated
by the strong solutions.

Finally, some numerical simulations are performed.

\end{abstract}



\bigskip

\section{Introduction}

The three-dimensional Navier-Stokes equations with damping have been studied
intensively over the last years. They describe the situation where there
exists resistance to the motion of a flow. One outstanding model in which a
damping term appears comes from the flow of cerebrospinal fluid inside the
porous brain tissues \cite{Linninger}. Such dissipative damping is also common
in many different models. For example, compressible Euler equations with
damping describe the flow of a compressible gas through a porous medium
\cite{HuangPan}, whereas Saint-Venant equations are used in oceanography to
describe the flow of viscous shallow water with friction \cite{Bresch}.

In this paper we study the asymptotic behaviour of weak solutions of the
following equation%
\begin{equation}
u_{t}-\mu\Delta u+\left(  u\text{\textperiodcentered}\nabla\right)
u+\alpha\left\vert u\right\vert ^{\beta-1}u+\nabla p=f,\text{ }\left(
x,t\right)  \in\Omega\times\left(  0,T\right)  , \label{Eq}%
\end{equation}
where $\Omega\subset\mathbb{R}^{3}$, $\beta\geq1,\ \mu,\alpha>0$, $\mu>0$ is
the kinematic viscosity and $u$ is the velocity vector of an incompressible
fluid satisfying Dirichlet boundary conditions.

We would like to point out that the damping term is very helpful from the
mathematical point of view, as it allows us to obtain solutions more regular
than in the standard Navier-Stokes equations without damping (that is, when
$\alpha=0$). For this reason it is possible to prove the existence of global
attractors for (\ref{Eq}), at least in a given range of the parameter $\beta$,
whereas to date this problem remains open for the standard three-dimensional
Navier-Stokes equations.

Existence of weak solutions for problem (\ref{Eq}) with initial condition in
the space of free-divergence square integrable functions was established at
first in \cite[Theorem 1]{CaiJiu} for $\beta\geq1$ and $\Omega=\mathbb{R}^{3}$
and in \cite[Theorem 2.1]{SongHou} for bounded domains $\Omega$. The
uniqueness of weak solutions was established in \cite{PardoValGim} for
$\beta\geq4$ and this result was extended in \cite{KimLiKim} for $\beta>3$ and
$\beta=3,\ \alpha\mu\geq\frac{1}{4}.$ A conditional result about smoothness of
weak solutions is given in \cite{WangZhou}.

Concerning global strong solutions of (\ref{Eq}) with $\Omega$ bounded and
more regular initial conditions, its existence was established in
\cite[Theorem 1.1]{KimLi} for either $\beta\in\left(  3,5\right)  $ or
$\beta=3,\ \ \alpha\mu>\frac{1}{4}.$ Recently, this result has been extended
to $\beta=5$ \cite{LiKimKimO}. Previously, the existence of strong solutions
was stated in \cite{SongHou} and \cite{SongLiangWu} for $\beta>3$. However, as
pointed out in \cite{KimLi}, \cite{LiKimKimO} it is unclear whether the proof
of this result is correct, because the function $-\Delta u$ is used as a test
function, which it seems cannot be done when $\Omega\not =\mathbb{R}^{3}$.
Instead, the function $Au$, where $A$ is the Stokes operator, has to be used
(see \cite[Theorem 1.1]{KimLi}).

Concerning global strong solutions of (\ref{Eq}) with $\Omega=\mathbb{R}^{3}$
and more regular initial conditions, its existence was established in
\cite{CaiJiu}, \cite{ZhangWuLu} for $\beta>3$, and in \cite{Zhou} for
$\beta=3$, $\alpha=\mu=1$. Also, in \cite{Zhou} uniqueness of strong solutions
was proved to be true for all $\beta\geq1$. If either the initial condition
$u_{0}$ is small enough or the viscosity $\mu$ is large, existence of global
strong solutions was stated in \cite{Zhong} for $1\leq\beta<3.$ It is pointed
out in \cite{LiKimKimO} that these results on existence of strong solutions
are conditional in the sense that it is necessary to prove previously the
existence of local strong solutions, because Galerkin approximations cannot be
used when $\Omega=\mathbb{R}^{3}$. This is done in \cite{LiKimKimO}.

The asymptotic behaviour of strong solutions when $3<\beta\leq5$ was studied
in \cite{SongHou}, \cite{SongHou2} and \cite{SongLiangWu} in the autonomous
and nonautonomous situation, proving the existence of the global attractor.
Again, the proof of these results is at least unclear as the test function
$-\Delta u$ has been used to obtain the suitable estimates of the solutions.
In this paper, we give an alternative proof of the existence of the strong
global attractor for either $3<\beta<5$ or $\beta=3,\ \alpha\mu>\frac{1}{4}$.
For $\beta\geq5$ the problem remains open.

For periodic boundary conditions the existence of weak and strong solutions
and global attractors has been studied in \cite{GautamMohan}, \cite{HajduRob},
\cite{KinraMoan}, \cite{MarTiTra}. 

The present paper is an improvement of \cite{PardoValGim}. In
\cite{PardoValGim} the existence of the global attractor when $\beta\geq3$ was
proved for weak solutions, which was the first result of this kind in the
literature. It is worth mentioning that it was necesary to use the theory of
multivalued semiflows for some values of the parameter $\beta$ due to the
absense of uniqueness. On top of that, the solutions and the attractor were
shown to be more regular for $\beta\geq4$. However, the proof uses the
function $-\Delta u$ a a test function, so as said before the proof is
unclear. In this paper, we prove the regularity results for weak solutions and
the global attractor for either $3<\beta<5$ or $\beta=3,\ \alpha\mu>\frac
{1}{4}.$ For $\beta\geq5$ the problem remains open.

This paper is organized as follows. In Section 2 we prove first suitable
estimates for weak solutions. In Section 3 we prove, for $\beta\geq3$, the
existence of the global attractor for the multivalued semiflow generated by
the weak solutions. When either $\beta>3$ or $\beta=3,\ \alpha\mu\geq\frac
{1}{4}$, this semiflow is a semigroup of operators and the attractor is proved
to be connected. If either $3<\beta<5$ or $\beta=3,\ \alpha\mu>\frac{1}{4}$,
we obtain more regularity of the attractor. Finally, when either $3<\beta<5$
or $\beta=3,\ \alpha\mu>\frac{1}{4}$, we establish that the global attractor
for the weak solutions is also the global attractor for the semigroupo
generated by the strong solutions.

\section{Estimates of weak solutions}

Consider a bounded open set $\Omega\subset\mathbb{R}^{3}$ with smooth boundary
$\partial\Omega$. We study the three-dimensional Navier-Stokes equations with
damping%
\begin{equation}
\left\{
\begin{array}
[c]{c}%
u_{t}-\mu\Delta u+\left(  u\text{\textperiodcentered}\nabla\right)
u+\alpha\left\vert u\right\vert ^{\beta-1}u+\nabla p=f,\text{ }\left(
x,t\right)  \in\Omega\times\left(  0,T\right)  ,\\
\operatorname{div}u=0,\ \left(  x,t\right)  \in\Omega\times\left(  0,T\right)
,\\
u\mid_{\partial\Omega}=0,\ t\in\left(  0,T\right)  ,\\
u\mid_{t=0}=u_{0}\text{, }x\in\Omega,
\end{array}
\right.  \label{NS}%
\end{equation}
where $\mu>0$ is the kinematic viscosity and $f$ is an external force. Also,
$\beta\geq1$ and $\alpha>0$ are given constants. The functions $u(x,t)=(u_{1}%
(x,t),u_{2}(x,t),u_{3}(x,t))$, $p(x,t)$ stand for the velocity field and the
pressure, respectively. Here and further, $\left\vert
\text{\textperiodcentered}\right\vert $ denotes in general the norm in
$\mathbb{R}^{d}$ for any $d\geq1$.

We define the usual function spaces
\begin{align*}
\mathcal{V}  &  =\{u\in(C_{0}^{\infty}(\Omega))^{3}:div\ u=0\},\ \\
H  &  =cl_{(L^{2}(\Omega))^{3}}\mathcal{V},\ \\
V  &  =cl_{(H_{0}^{1}(\Omega))^{3}}\mathcal{V},\
\end{align*}
where $cl_{X}$ denotes the closure in the space $X$. It is well known that
$H,\ V$ are separable Hilbert spaces and identifying $H$ and its dual we have
$V\subset H\subset V^{\prime}$ with dense and continuous injections. We denote
by $(\cdot,\cdot),\ \left\vert \text{\textperiodcentered}\right\vert $ and
$((\cdot,\cdot)),\ \Vert\cdot\Vert$ the inner product and norm in $H$ and $V$,
respectively, and by $\left\langle \text{\textperiodcentered}%
,\text{\textperiodcentered}\right\rangle $ duality between $V^{\prime}$ and
$V$. Let $H_{w}$ be the space $H$ endowed with the weak topology. As usual, we
define the continuous trilinear form $b:V\times V\times V\rightarrow
\mathbb{R}$ by
\[
b(u,v,w)=\sum\limits_{i,j=1}^{3}\int\limits_{\Omega}u_{i}\frac{\partial v_{j}%
}{\partial x_{i}}w_{j}dx.
\]
It is well-known that $b(u,v,v)=0$, if $u\in V,v\in\left(  H_{0}^{1}\left(
\Omega\right)  \right)  ^{3}$. For $u,v\in V$ we denote by $B\left(
u,v\right)  $ the element of $V^{\prime}$ defined by $\left\langle B\left(
u,v\right)  ,w\right\rangle =b(u,v,w)$, for all $w\in V$.

The norm in the spaces $L^{p}\left(  \Omega\right)  $, $\left(  L^{p}\left(
\Omega\right)  \right)  ^{3}$, $p\geq1$, will be denoted indistinctly by
$\left\vert \text{\textperiodcentered}\right\vert _{p}.$

Let $P$ be the orthogonal projection from $(L^{2}(\Omega))^{3}$ onto $H$ and
$Au=-P\Delta u$ be the Stokes operator, defined by $\left\langle
Au,v\right\rangle =((u,v))$ for $u,v\in V.$ Since the boundary $\partial
\Omega$ is smooth, $D(A)=\left(  H^{2}(\Omega)\right)  ^{3}\cap V$ and
$\left\Vert Au\right\Vert _{2}$ defines a norm in $D(A)$ which is equivalent
to the norm in $\left(  H^{2}(\Omega)\right)  ^{3}$.

For $u_{0}\in H,\ f\in H$ the function
\begin{equation}
u\in L^{\infty}\left(  0,T;H\right)  \cap L^{2}\left(  0,T;V\right)  \cap
L^{\beta+1}\left(  0,T;\left(  L^{\beta+1}\left(  \Omega\right)  \right)
^{3}\right)  \label{RegSol}%
\end{equation}
is said to be a weak solution to problem (\ref{NS}) on $\left(  0,T\right)  $
if $u\left(  0\right)  =u_{0}$ and
\begin{equation}
\frac{d}{dt}\left(  u,v\right)  +\mu(\left(  u,v\right)  )+b\left(
u,u,v\right)  +\alpha\left(  \left\vert u\right\vert ^{\beta-1}u,v\right)
=\left(  f,v\right)  , \label{EqSol}%
\end{equation}
for any $v\in V\cap\left(  L^{\beta+1}\left(  \Omega\right)  \right)  ^{3}$,
in the sense of scalar distributions.

We recall the following well-known result on existence of weak solutions.

\begin{theorem}
\label{ExistWeakSol}\cite[Theorem 2.1]{SongHou} For any $u_{0}\in H,\ f\in H$,
$\beta\geq1$ there exists at least one weak solution $u$ to problem (\ref{NS}).
\end{theorem}

Let $Y=V^{\prime}+\left(  L^{\frac{\beta+1}{\beta}}\left(  \Omega\right)
\right)  ^{3}$, the dual space of $V\cap\left(  L^{\beta+1}\left(
\Omega\right)  \right)  ^{3}$. We note that by standard estimates on $B$ for
any weak solution we have that%
\begin{align*}
Au  &  \in L^{2}(0,T;V^{\prime}),\\
B(u,u)  &  \in L^{\frac{4}{3}}(0,T;V^{\prime}),\\
\left\vert u\right\vert ^{\beta-1}u  &  \in L^{\frac{\beta+1}{\beta}}\left(
0,T;\left(  L^{\beta+1}\left(  \Omega\right)  \right)  ^{3}\right)  ,
\end{align*}
which implies in particular that
\begin{align*}
-\mu Au-B(u,u)-\alpha\left\vert u\right\vert ^{\beta-1}u+f  &  \in L^{\frac
{4}{3}}\left(  0,T;V^{\prime}\right)  +L^{\frac{\beta+1}{\beta}}\left(
0,T;\left(  L^{\frac{\beta+1}{\beta}}\left(  \Omega\right)  \right)
^{3}\right) \\
&  \subset L^{1}(0,T;Y).
\end{align*}
It follows from equality (\ref{EqSol}) and a standard result \cite[p.250,
Lemma 1.1]{Temam} that%
\begin{equation}
\frac{du}{dt}=-\mu Au-B(u,u)-\alpha\left\vert u\right\vert ^{\beta-1}u+f
\label{EqSol2}%
\end{equation}
in the sense of $Y$-valued distributions. Hence, the derivate $\dfrac{du}{dt}$
belongs to the space%
\[
L^{\frac{4}{3}}\left(  0,T;V^{\prime}\right)  +L^{\frac{\beta+1}{\beta}%
}\left(  0,T;\left(  L^{\frac{\beta+1}{\beta}}\left(  \Omega\right)  \right)
^{3}\right)
\]
and equality (\ref{EqSol2}) is satisfied in the space $Y$ for a.a.
$t\in\left(  0,T\right)  .$

In order to obtain good estimates of weak solutions we need $\dfrac{du}{dt}$
to be more regular. We can obtain such a result for $\beta\geq3$.

\begin{lemma}
\label{PropertiesWeakSol}Let $u$ be a weak solution to (\ref{NS}) such that
$u\in L^{q}\left(  0,T;\left(  L^{q}\left(  \Omega\right)  \right)
^{3}\right)  $ with $q\geq4$. Then%
\begin{equation}
\frac{du}{dt}\in L^{2}\left(  0,T;V^{\prime}\right)  +L^{\frac{\beta+1}{\beta
}}\left(  0,T;\left(  L^{\frac{\beta+1}{\beta}}\left(  \Omega\right)  \right)
^{3}\right)  , \label{PropDeriv}%
\end{equation}%
\begin{equation}
u\in C\left(  [0,T],H\right)  , \label{ContWeakSol}%
\end{equation}
the map $t\mapsto\left\Vert u\left(  t\right)  \right\Vert _{H}^{2}$ is
absolutely continuous and%
\begin{equation}
\frac{d}{dt}\left\vert u\left(  t\right)  \right\vert ^{2}=2\left\langle
u,\frac{du}{dt}\right\rangle \text{ for a.a. }t\in\left(  0,T\right)  .
\label{NormDeriv}%
\end{equation}

\end{lemma}

\begin{proof}
Using the well-known inequality (see \cite[p.297]{Temam})
\[
\left\vert b\left(  u,u,v\right)  \right\vert \leq C\left\vert u\right\vert
_{4}^{2}\left\Vert v\right\Vert \text{, }\forall u,v\in V,
\]
and $u\in L^{q}\left(  0,T;\left(  L^{q}\left(  \Omega\right)  \right)
^{3}\right)  \subset L^{4}\left(  0,T;\left(  L^{4}\left(  \Omega\right)
\right)  ^{3}\right)  $, we have%
\[
B\left(  u,u\right)  \in L^{2}\left(  0,T;V^{\prime}\right)  ,
\]
so (\ref{PropDeriv}) follows.

Properties (\ref{ContWeakSol})-(\ref{NormDeriv}) follow from \cite[Chapter II,
Theorem 1.8]{ChepVishikBook}.
\end{proof}

\begin{corollary}
If $\beta\geq3$, then any weak solution to (\ref{NS}) satisfies
(\ref{PropDeriv})-(\ref{NormDeriv}).
\end{corollary}

\begin{proof}
Since $u\in L^{\beta+1}\left(  0,T;\left(  L^{\beta+1}\left(  \Omega\right)
\right)  ^{3}\right)  $ and $\beta\geq3$, we obtain that $u$ belongs to
$L^{q}\left(  0,T;\left(  L^{q}\left(  \Omega\right)  \right)  ^{3}\right)  $
with $q\geq4$.
\end{proof}


\begin{lemma}
\label{EstWeakSol}If $\beta\geq3$, then any weak solution satisfies the
estimates%
\begin{equation}
\left\vert u(t)\right\vert ^{2}\leq e^{-\mu\lambda_{1}t}\left\vert
u_{0}\right\vert ^{2}+\frac{\left\vert f\right\vert ^{2}}{\mu^{2}\lambda
_{1}^{2}}, \label{EstWeak1}%
\end{equation}%
\begin{equation}
\mu\int_{s}^{t}\left\Vert u\right\Vert ^{2}d\tau+2\alpha\int_{s}^{t}\left\vert
u\right\vert _{\beta+1}^{\beta+1}d\tau\leq\left\vert u_{0}\right\vert
^{2}+\frac{\left\vert f\right\vert ^{2}}{\mu^{2}\lambda_{1}^{2}}+\frac{1}%
{\mu\lambda_{1}}\left\vert f\right\vert ^{2}\left(  t-s\right)  ,
\label{EstWeak2}%
\end{equation}
for any $t\geq s\geq0.$
\end{lemma}

\begin{proof}
Multiplying equality (\ref{EqSol2}) by $u$ and using (\ref{NormDeriv}) and
$b\left(  u,u,u\right)  =0$ we have%
\begin{equation}
\frac{1}{2}\frac{d}{dt}\left\vert u\right\vert ^{2}+\mu\left\Vert u\right\Vert
^{2}+\alpha\left\vert u\right\vert _{\beta+1}^{\beta+1}=\left(  f,u\right)
\leq\frac{\mu\lambda_{1}}{2}\left\vert u\right\vert ^{2}+\frac{1}{2\mu
\lambda_{1}}\left\vert f\right\vert ^{2}. \label{Ineq1}%
\end{equation}
As $\mu\left\Vert u\right\Vert ^{2}\geq\mu\lambda_{1}\left\vert u\right\vert
^{2}$, we deduce that%
\begin{equation}
\frac{d}{dt}\left\vert u\right\vert ^{2}+\mu\lambda_{1}\left\vert u\right\vert
^{2}\leq\frac{1}{\mu\lambda_{1}}\left\vert f\right\vert ^{2} \label{Ineq2}%
\end{equation}
and Gronwall's lemma yields%
\[
\left\vert u(t)\right\vert ^{2}\leq e^{-\mu\lambda_{1}t}\left\vert
u_{0}\right\vert ^{2}+\frac{\left\vert f\right\vert ^{2}}{\mu^{2}\lambda
_{1}^{2}}.
\]
By $\mu\left\Vert u\right\Vert ^{2}\geq\frac{\mu}{2}\left\Vert u\right\Vert
^{2}+\frac{\mu\lambda_{1}}{2}\left\vert u\right\vert ^{2}$, integrating over
the interval $\left(  s,t\right)  $ in (\ref{Ineq1}) it follows that%
\begin{align*}
\mu\int_{s}^{t}\left\Vert u\right\Vert ^{2}d\tau+2\alpha\int_{s}^{t}\left\vert
u\right\vert _{\beta+1}^{\beta+1}d\tau &  \leq\left\vert u\left(  s\right)
\right\vert ^{2}+\frac{1}{\mu\lambda_{1}}\left\vert f\right\vert ^{2}\left(
t-s\right) \\
&  \leq e^{-\mu\lambda_{1}s}\left\vert u_{0}\right\vert ^{2}+\frac{\left\vert
f\right\vert ^{2}}{\mu^{2}\lambda_{1}^{2}}+\frac{1}{\mu\lambda_{1}}\left\vert
f\right\vert ^{2}\left(  t-s\right)  ,
\end{align*}
so the lemma is proved.
\end{proof}

\bigskip

The uniqueness of weak solution was established at first in \cite[Theorem
2.4]{PardoValGim} for $\beta\geq4.$ Later on, in \cite[Corollary
2.1]{KimLiKim} this result was extended for $\beta>3$ and $\beta
=3,\ 4\alpha\mu\geq1.$

\begin{theorem}
\label{UniquenessWeak}Let either $\beta>3$ or $\beta=3$ and $4\alpha\mu\geq1$.
Then for any $u_{0}\in H$ there exists a unique weak solution $u\left(
\text{\textperiodcentered}\right)  $ to problem (\ref{NS}), which is
continuous with respect to the initial datum $u_{0}$.
\end{theorem}

\section{Global attractor for weak and strong solutions}

Our aim now is to prove the existence of the global attractor for the weak and
strong solutions of problem (\ref{NS}).

We shall divide this section into two cases: 1) $\beta\geq3$; 2) $\beta>3$ or
$\beta=3,\ 4\alpha\mu\geq1$. In the first one, as in general more than one
solution can possibly exist for a given initial datum, we make use of the
theory of attractors for multivalued semiflows to prove the existence of a
global attractor. In the second one, uniqueness of weak solutions implies that
we can define a semigroup of operators, to which we can apply the classical
theory of attractors for semigroups, proving the existence of a global
connected attractor. More regularity of the attractor is obtained if either
$\beta>3$ or $\beta=3,\ 4\alpha\mu>1$. The attractor is shown to be the global
attractor for the strong solutions as well.

\subsection{Case 1: $\beta\geq3$}

Let us define the set
\[
D_{T}\left(  u_{0}\right)  =\{u\left(  \text{\textperiodcentered}\right)
\text{ is a weak solution of (\ref{NS}) in the interval }\left(  0,T\right)
\}.
\]
We know by Theorem \ref{ExistWeakSol} that for any $u_{0}\in H$ and $T>0$ the
set $D_{T}\left(  u_{0}\right)  $ is non-empty.

We observe that as $q=\frac{\beta+1}{\beta}\leq\frac{4}{3}$, we have the time
derivative of a weak solution satisfies
\[
\frac{du}{dt}\in L^{\frac{4}{3}}\left(  0,T;V^{\prime}\right)  +L^{q}%
(0,T;\left(  L^{q}(\Omega)\right)  ^{3})\subset L^{q}(0,T;V^{\prime}+\left(
L^{q}(\Omega)\right)  ^{3}).
\]

\begin{lemma}
\label{Concatenation}Let $\beta\geq3$. If $u\left(  \text{\textperiodcentered
}\right)  \in D_{T}\left(  u_{0}\right)  $, then for any $s\in\left(
0,T\right)  $, the function $w\left(  \text{\textperiodcentered}\right)
=u\left(  \text{\textperiodcentered}+s\right)  $ belongs to $D_{T-s}\left(
u\left(  s\right)  \right)  .$

If $u\left(  \text{\textperiodcentered}\right)  \in D_{s}\left(  u_{0}\right)
$ and $w\left(  \text{\textperiodcentered}\right)  \in D_{T-s}\left(  u\left(
s\right)  \right)  $, then the function%
\[
z\left(  t\right)  =\left\{
\begin{array}
[c]{c}%
u\left(  t\right)  \text{ if }t\in\lbrack0,s],\\
w\left(  t-s\right)  \text{ if }t\in\lbrack s,T],
\end{array}
\right.
\]
belongs to $D_{T}\left(  u_{0}\right)  .$
\end{lemma}

\begin{proof}
Let $u\left(  \text{\textperiodcentered}\right)  \in D_{T}\left(
u_{0}\right)  $. Then it is obvious that
\begin{equation}
w\left(  \text{\textperiodcentered}\right)  =u\left(
\text{\textperiodcentered}+s\right)  \in L^{\infty}\left(  0,T-s;H\right)
\cap L^{2}\left(  0,T-s;V\right)  \cap L^{\beta+1}\left(  0,T-s;L^{\beta
+1}\left(  \Omega\right)  \right)  . \label{Prop1Sol}%
\end{equation}
Also, (\ref{EqSol}) implies that for any $v\in V\cap\left(  L^{\beta+1}\left(
\Omega\right)  \right)  ^{3}$, $\phi\in C_{0}^{\infty}(0,T-s)$ one has%
\begin{align*}
&  -\int_{0}^{T-s}\left(  w(\tau),v\right)  \phi^{\prime}(\tau)d\tau\\
&  +\int_{0}^{T-s}\mu(\left(  w(\tau),v\right)  )+b\left(  w(\tau
),w(\tau),v\right)  +\alpha\left(  \left\vert w(\tau)\right\vert ^{\beta
-1}w(\tau),v\right)  \phi(\tau)d\tau\\
&  =-\int_{s}^{T}\left(  u(r),v\right)  \phi^{\prime}(r-s)dr\\
&  +\int_{s}^{T}\mu(\left(  u(r),v\right)  )+b\left(  u(r),u(r),v\right)
+\alpha\left(  \left\vert u(r)\right\vert ^{\beta-1}u(r),v\right)
\phi(r-s)dr\\
&  =\int_{s}^{T}\left(  f,v\right)  \phi(r-s)dr=\int_{0}^{T-s}\left(
f,v\right)  \phi(\tau)d\tau,
\end{align*}
so $w$ satisfies (\ref{EqSol}) in the interval $\left(  0,T-s\right)  .$ We
infer that $w\in D_{T-s}\left(  u\left(  s\right)  \right)  .$

Let now $u\left(  \text{\textperiodcentered}\right)  \in D_{s}\left(
u_{0}\right)  $ and $w\left(  \text{\textperiodcentered}\right)  \in
D_{T-s}\left(  u\left(  s\right)  \right)  $. Arguing as in the previous case
we obtain that $w(t-s)$ satisfies equality (\ref{EqSol}) in the interval
$\left(  s,T\right)  $. As the time derivative of a weak solution belongs to
$L^{q}\left(  0,T;V^{\prime}+\left(  L^{q}(\Omega)\right)  ^{3}\right)  $, by
\cite[p.250, Lemma 1.1]{Temam} equality (\ref{EqSol}) is equivalent to saying
that%
\[
\int_{0}^{T}\left(  \left\langle \frac{du}{dt},\xi\right\rangle +\mu\left(
\left(  u,\xi\right)  \right)  +\left\langle B(u,u),\xi\right\rangle \right)
dt+\alpha\int_{0}^{T}\int_{\Omega}\left\vert u\right\vert ^{\beta-1}u\xi
dxdt=\int_{0}^{T}\left(  f,\xi\right)  dt,
\]
for any $\xi\in L^{\beta+1}\left(  0,T;V\cap\left(  L^{\beta+1}\left(
\Omega\right)  \right)  ^{3}\right)  $. The function $z$ satisfies
(\ref{RegSol}) in the interval $\left(  0,T\right)  $ and this equality as
well. Indeed, denoting $h(t)=w(t-s)$ we have%
\begin{align*}
&  \int_{0}^{T}\left(  \left\langle \frac{dz}{dt},\xi\right\rangle +\mu\left(
\left(  z,\xi\right)  \right)  +\left\langle B(z,z),\xi\right\rangle \right)
dt+\alpha\int_{0}^{T}\int_{\Omega}\left\vert z\right\vert ^{\beta-1}z\xi
dxdt\\
&  =\int_{0}^{s}\left(  \left\langle \frac{du}{dt},\xi\right\rangle
+\mu\left(  \left(  u,\xi\right)  \right)  +\left\langle B(u,u),\xi
\right\rangle \right)  dt+\alpha\int_{0}^{s}\int_{\Omega}\left\vert
u\right\vert ^{\beta-1}u\xi dxdt\\
&  +\int_{s}^{T}\left(  \left\langle \frac{dh}{dt},\xi\right\rangle
+\mu\left(  \left(  h,\xi\right)  \right)  +\left\langle B(h,h),\xi
\right\rangle \right)  dt+\alpha\int_{s}^{T}\int_{\Omega}\left\vert
h\right\vert ^{\beta-1}h\xi dxdt\\
&  =\int_{0}^{s}\left(  f,\xi\right)  dt+\int_{s}^{T}\left(  f,\xi\right)
dt=\int_{0}^{T}\left(  f,\xi\right)  dt,
\end{align*}
proving that $z$ is really a weak solution.
\end{proof}

\bigskip

In view of this lemma every solution can be extended to a globally defined
one, that is, a solution which exists for $t\in\lbrack0,+\infty)$. In this
situation we denote by $D\left(  u_{0}\right)  $ the set of all globally
defined solutions with initial condition $u_{0}$ and observe that for any
$t\geq0$ the following equality holds:%
\[
\{u(t):u\in D\left(  u_{0}\right)  \}=\{u(t):u\in\cup_{T>0}D_{T}\left(
u_{0}\right)  \}.
\]

Denote by $P\left(  H\right)  $ the set of all non-empty subsets of $H$. Let
us define the following (possibly multivalued) family of operators
$G:\mathbb{R}^{+}\times H\rightarrow P\left(  H\right)  :$%
\[
G\left(  t,u_{0}\right)  =\{y\in H:y=u\left(  t\right)  \text{, }u\left(
\text{\textperiodcentered}\right)  \in D\left(  u_{0}\right)  \}.
\]
Using Lemma \ref{Concatenation} we can easily prove that $G$ is a strict
multivalued semiflow, that is, the following two properties hold:

\begin{itemize}
\item $G\left(  0,u_{0}\right)  =u_{0}$ for all $u_{0}\in H;$

\item $G\left(  t+s,u_{0}\right)  =G\left(  t,G\left(  s,u_{0}\right)
\right)  $, for all $u_{0}\in H$, $t,s\geq0.$
\end{itemize}

The set $\mathcal{A}$ is a global attractor for $G$ if:

\begin{itemize}
\item $\mathcal{A}$ is negatively invariant, i.e., $\mathcal{A}\subset
G\left(  t,\mathcal{A}\right)  $ for all $t\geq0;$

\item $\mathcal{A}$ attracts every bounded set of $H$, that is,%
\[
dist\left(  G\left(  t,B\right)  ,\mathcal{A}\right)  \rightarrow0\text{ as
}t\rightarrow+\infty.
\]

\end{itemize}

It is invariant if, moreover, $\mathcal{A}=G\left(  t,\mathcal{A}\right)  $
for all $t\geq0.$

The next lemma is crucial for proving the existence of a global attractor.

\begin{lemma}
\label{LemmaConverg}Assume that $\beta\geq3$. Let $u_{0}^{n}\rightarrow u_{0}$
weakly in $H$ and let $u_{n}\left(  \text{\textperiodcentered}\right)  \in
D\left(  u_{0}^{n}\right)  $. Then there exists a weak solution $u\left(
\text{\textperiodcentered}\right)  $ to (\ref{NS}) with $u\left(  0\right)
=u_{0}$ and a subsequence $u_{n_{k}}\left(  \text{\textperiodcentered}\right)
$ such that $u_{n_{k}}\rightarrow u$ in $C\left(  [\varepsilon,T],H\right)  $
for all $0<\varepsilon<T.$

If, moreover, $u_{0}^{n}\rightarrow u_{0}$ strongly in $H$, then $u_{n_{k}%
}\rightarrow u$ in $C\left(  [0,T],H\right)  $ for all $T>0.$
\end{lemma}

\begin{proof}
We fix $T>0$. We deduce from Lemma \ref{EstWeakSol} that the sequence $u_{n}$
is bounded in
\[
L^{\infty}\left(  0,T;H\right)  \cap L^{2}\left(  0,T;V\right)  \cap
L^{\beta+1}\left(  0,T;\left(  L^{\beta+1}\left(  \Omega\right)  \right)
^{3}\right)  .
\]
Also, using (\ref{EqSol2}) and standard estimates (see \cite[p.297]{Temam}) we
obtain that $\dfrac{du_{n}}{dt}$ is bounded in the space $L^{q}\left(
0,T;V^{\prime}+\left(  L^{q}(\Omega)\right)  ^{3}\right)  $.

Thus, making use of the compactness theorem \cite{Lions} we obtain a function
$u\left(  \text{\textperiodcentered}\right)  $ and a subsequence (denoted
again by $u_{n}$) such that%
\begin{align}
u_{n}  &  \rightarrow u\text{ weakly star in }L^{\infty}\left(  0,T;H\right)
,\label{Converg}\\
u_{n}  &  \rightarrow u\text{ weakly in }L^{2}\left(  0,T;V\right)
,\nonumber\\
u_{n}  &  \rightarrow u\text{ weakly in }L^{\beta+1}\left(  0,T;\left(
L^{\beta+1}\left(  \Omega\right)  \right)  ^{3}\right)  ,\nonumber\\
\frac{du_{n}}{dt}  &  \rightarrow\frac{du}{dt}\text{ weakly in }L^{q}\left(
0,T;V^{\prime}+\left(  L^{q}(\Omega)\right)  ^{3}\right)  ,\nonumber\\
u_{n}  &  \rightarrow u\text{ strongly in }L^{2}\left(  0,T;H\right)
,\nonumber\\
u_{n}\left(  t,x\right)   &  \rightarrow u\left(  t,x\right)  \text{ for a.a.
}\left(  t,x\right)  .\nonumber
\end{align}
Let us prove that
\begin{equation}
u_{n}\left(  t_{n}\right)  \rightarrow u\left(  t_{0}\right)  \text{ weakly in
}H \label{ConvergWeak}%
\end{equation}
for any sequence $\{t_{n}\}$ such that $t_{n}\rightarrow t_{0}$, where
$t_{n},t_{0}\in\lbrack0,T]$. The time derivatives are bounded in the space
$L^{q}\left(  0,T;V^{\prime}+\left(  L^{q}(\Omega)\right)  ^{3}\right)  $,
which implies readily that the sequence $u_{n}\left(
\text{\textperiodcentered}\right)  $ is equicontinuous in the space
$V^{\prime}+\left(  L^{q}(\Omega)\right)  ^{3}$. Moreover, $u_{n}\left(
t_{n}\right)  $ is bounded in $H$, and then the compact embedding $H\subset
V^{\prime}$ yields that it is relatively compact in $V^{\prime}+\left(
L^{q}(\Omega)\right)  ^{3}$. Hence, by Ascoli-Arzel\`{a}'s theorem we have
$u_{n}\rightarrow u$ in $C\left(  [0,T],V^{\prime}+\left(  L^{q}%
(\Omega)\right)  ^{3}\right)  $. Thus, by a contradiction argument we obtain
that $u_{n}\left(  t_{n}\right)  \rightarrow u\left(  t_{0}\right)  $ weakly
in $H$. In particular, we have that $u\left(  0\right)  =u_{0}$.

Further, we need to check that $u\left(  \text{\textperiodcentered}\right)  $
is a weak solution to problem (\ref{NS}).

The sequence $h\left(  u_{n}\left(  \text{\textperiodcentered}\right)
\right)  =\left\vert u_{n}\left(  \text{\textperiodcentered}\right)
\right\vert ^{\beta-1}u_{n}\left(  \text{\textperiodcentered}\right)  $ is
bounded in $L^{q}\left(  0,T;\left(  L^{q}\left(  \Omega\right)  \right)
^{3}\right)  $ and $h\left(  u_{n}\left(  t,x\right)  \right)  \rightarrow
h\left(  u\left(  t,x\right)  \right)  $ for a.a. $\left(  t,x\right)  $.
Hence, $h\left(  u_{n}\left(  \text{\textperiodcentered}\right)  \right)
\rightarrow h\left(  u\left(  \text{\textperiodcentered}\right)  \right)  $
weakly in $L^{q}\left(  0,T;\left(  L^{q}\left(  \Omega\right)  \right)
^{3}\right)  $ \cite[Lemma 8.3]{Robinson}.

In order to show that $u$ is a weak solution it remains to pass to the limit
in the term $B$. Since $u_{n}\rightarrow u$ in $L^{2}(0,T;H)$ implies that
$u_{ni}u_{nj}\rightarrow u_{i}u_{j}$ in $L^{1}(0,T;L^{1}(\Omega))$, for any
$\zeta\in\mathcal{V}$, $\phi\in C_{0}^{\infty}(0,T)$ we have%
\begin{align*}
\int_{0}^{T}\left(  b(u_{n},u_{n},\zeta)-b(u,u,\zeta)\right)  \phi dt  &
=-\int_{0}^{T}\left(  b(u_{n},\zeta,u_{n})-b(u,\zeta,u)\right)  \phi dt\\
&  =-\sum_{i,j=1}^{3}\int_{0}^{T}\int_{\Omega}\left(  u_{ni}u_{nj}-u_{i}%
u_{j}\right)  \frac{\partial\zeta_{j}}{\partial x_{i}}\phi dxdt\rightarrow0,
\end{align*}
as $n\rightarrow\infty$.

We conclude that equality (\ref{EqSol}) is satisfied for the function $u$ for
all $\zeta\in\mathcal{V}$, and by density of $\mathcal{V}$ in $V$ we obtain
that (\ref{EqSol}) holds true. Thus, $u$ is a weak solution.

Finally, we will prove that $u_{n}\rightarrow u$ in $C\left(  [\varepsilon
,T],H\right)  $ for all $0<\varepsilon<T$. From (\ref{Ineq2}) we get
\[
\left\vert u_{n}\left(  t\right)  \right\vert ^{2}\leq\left\vert u_{n}\left(
s\right)  \right\vert ^{2}+\frac{1}{\mu\lambda_{1}}\left\vert f\right\vert
^{2}\left(  t-s\right)  \text{, for any }s\leq t,
\]
and the same inequality is true for $u$. Hence, the functions $J_{n}\left(
t\right)  =\left\vert u_{n}\left(  t\right)  \right\vert ^{2}-\frac{1}%
{\mu\lambda_{1}}\left\vert f\right\vert ^{2}t,\ J\left(  t\right)  =\left\vert
u\left(  t\right)  \right\vert ^{2}-\frac{1}{\mu\lambda_{1}}\left\vert
f\right\vert ^{2}t$ are non-increasing and continuous. Take a sequence
$t_{n}\rightarrow t_{0}$ with $t_{n},t_{0}\in\lbrack\varepsilon,T]$. We know
that $u_{n}\left(  t_{n}\right)  \rightarrow u\left(  t_{0}\right)  $ weakly
in $H$, so
\begin{equation}
\left\vert u\left(  t_{0}\right)  \right\vert \leq\lim\inf\ \left\vert
u\left(  t_{n}\right)  \right\vert . \label{LimInf}%
\end{equation}
\ It is a consequence of (\ref{Converg}) that $J_{n}\left(  t\right)
\rightarrow J\left(  t\right)  $ for a.a. $t$. Then we can choose $t_{k}%
<t_{0}$ as close to $t_{0}$ as we wish such that $J_{n}\left(  t_{k}\right)
\rightarrow J\left(  t_{k}\right)  $, and we can assume without loss of
generality that $t_{k}<t_{n}$. Therefore,%
\begin{align*}
J_{n}\left(  t_{n}\right)  -J\left(  t_{0}\right)   &  =J_{n}\left(
t_{n}\right)  -J_{n}\left(  t_{k}\right)  +J_{n}\left(  t_{k}\right)
-J\left(  t_{k}\right)  +J\left(  t_{k}\right)  -J\left(  t_{0}\right) \\
&  \leq\left\vert J_{n}\left(  t_{k}\right)  -J\left(  t_{k}\right)
\right\vert +\left\vert J\left(  t_{k}\right)  -J\left(  t_{0}\right)
\right\vert .
\end{align*}
Since $u\left(  \text{\textperiodcentered}\right)  $ is continuous, for any
$\delta>0$ there exists $t_{k}$ and $N\left(  t_{k}\right)  $ such that
$\left\vert J\left(  t_{k}\right)  -J\left(  t_{0}\right)  \right\vert
\leq\delta/2$ and $\left\vert J_{n}\left(  t_{k}\right)  -J\left(
t_{k}\right)  \right\vert \leq\delta/2$ for all $n\geq N$. This implies that%
\begin{equation}
\lim\sup\ \left\vert u\left(  t_{n}\right)  \right\vert \leq\left\vert
u\left(  t_{0}\right)  \right\vert . \label{LimSup}%
\end{equation}
Joining (\ref{LimInf}) and (\ref{LimSup}) we deduce that $\left\vert u\left(
t_{n}\right)  \right\vert \rightarrow\left\vert u\left(  t_{0}\right)
\right\vert $ and then $u\left(  t_{n}\right)  \rightarrow u\left(
t_{0}\right)  $ in $H.$

Since $T>0$ is arbitrary, by a diagonal arguments we obtain a common
subsequence on an arbitrary interval $[\varepsilon,T].$

The first part of the lemma is proved.

For the second part, we need to prove only that $u\left(  t_{n}\right)
\rightarrow u\left(  0\right)  $ if $t_{n}\rightarrow0$, $t_{n}\geq0$. For
this aim we repeat the above argument with $t_{k}=0=t_{0}$. Hence,%
\[
J_{n}\left(  t_{n}\right)  -J\left(  t_{0}\right)  =J_{n}\left(  t_{n}\right)
-J_{n}\left(  0\right)  +J_{n}\left(  0\right)  -J\left(  0\right)
\leq\left\vert J_{n}\left(  0\right)  -J\left(  0\right)  \right\vert
\rightarrow0\text{ as }n\rightarrow\infty,
\]
because $u_{n}\left(  0\right)  \rightarrow u\left(  0\right)  $ in $H$. Then
we obtain the result arguing in the same way as above.
\end{proof}


\begin{corollary}
\label{Compact}Assume that $\beta\geq3$. For any $t\geq0$ the map
$u_{0}\mapsto G(t,u_{0})$ has compact values and closed graph. In addition,
for any $t_{0}>0$ the map $G\left(  t_{0},\text{\textperiodcentered}\right)  $
is compact.
\end{corollary}

The map $u_{0}\mapsto G(t,u_{0})$ is said to be upper semicontinuous if for
all $u_{0}\in H$ and any neighborhood $O$ of $u_{0}$ in $H$ there exists
$\delta>0$ such that $G\left(  t,u\right)  \subset O$ for all $u$ satisfying
$\left\Vert u-u_{0}\right\Vert <\delta$.

\begin{lemma}
\label{USC}Assume that $\beta\geq3$. For any $t\geq0$ the map $u_{0}\mapsto
G(t,u_{0})$ is upper semicontinuous.
\end{lemma}

\begin{proof}
If not, there exist $u_{0}\in H,\ t>0$, sequences $u_{n}^{0}\rightarrow u_{0}%
$, $y_{n}\in G\left(  t,u_{0}^{n}\right)  $ and a neighborhood $O$ of
$G\left(  t,u_{0}\right)  $ such that $y_{n}\not \in O$. Let $y_{n}%
=u_{n}\left(  t\right)  $, where $u_{n}\left(  \text{\textperiodcentered
}\right)  \in D\left(  u_{0}^{n}\right)  $. Then by Lemma \ref{LemmaConverg}
there is a subsequence $y_{n_{k}}$ satisfying $y_{n_{k}}\rightarrow y\in
G\left(  t,u_{0}\right)  $, which is a contradiction.
\end{proof}

\bigskip

A set $B_{0}$ is called absorbing if for any bounded set $B$ there exists a
time $T\left(  B\right)  $ such that%
\[
G\left(  t,B\right)  \subset B_{0}\text{ for any }t\geq T.
\]
The semiflow $G$ is said to be asymptotically compact if for any bounded
subset $B$ every sequence $y_{k}\in G\left(  t_{k},B\right)  $, where
$t_{k}\rightarrow+\infty$, is relatively compact in $H.$

The following conditions are sufficient in order to obtain a global compact
invariant minimal attractor $\mathcal{A}$ for a strict multivalued semiflow
$G$ \cite[Theorem 3\ and Remark 8]{MelnikValero98}:

\begin{enumerate}
\item $G$ possesses a bounded absorbing set $B_{0};$

\item $G$ is asymptotically compact;

\item $G$ has closed values;

\item the map $u_{0}\mapsto G(t,u_{0})$ is upper semicontinuous.
\end{enumerate}

\begin{theorem}
\label{AttrExist}Assume that $\beta\geq3$. Then $G$ has a global invariant
compact attractor $\mathcal{A}$, which is minimal among all closed attracting sets.
\end{theorem}

\begin{proof}
We need to check the four aforementioned conditions.

It follows from (\ref{EstWeak1}) that the ball $B_{0}=\{u\in H:\left\Vert
u\right\Vert _{H}^{2}\leq1+\frac{\left\Vert f\right\Vert _{H}^{2}}{\mu
^{2}\lambda_{1}^{2}}\}$ is absorbing.

In view of Corollary \ref{Compact} and Lemma \ref{USC} $G$ has compact (and
then closed) values and the map $u_{0}\mapsto G(t,u_{0})$ is upper semicontinuous.

Finally, again by Corollary \ref{Compact} the operator $G\left(
1,\text{\textperiodcentered}\right)  $ is compact. Hence, for any bounded set
$B$ an arbitrary sequence $y_{n}\in G\left(  t_{n},B\right)  $, which belongs
to
\[
G\left(  1,G\left(  t_{n}-1,B\right)  \right)  \subset G\left(  1,B_{0}%
\right)  ,\ \text{for all } n\geq N,
\]
is relatively compact in $H$, so $G$ is asymptotically compact.
\end{proof}

\bigskip

We can give also some information about the structure of the global attractor
in terms of bounded complete trajectories, which are continuous functions
$\gamma:\mathbb{R}\rightarrow H$ such that $u\left(  \text{\textperiodcentered
}\right)  =\gamma\left(  \text{\textperiodcentered}+s\right)  $ belongs to
$D\left(  u_{0}\right)  $ for all $s\in\mathbb{R}$ and satisfying that
$\cup_{t\in\mathbb{R}}\gamma\left(  t\right)  $ is a bounded set. Indeed, by
Lemmas \ref{Concatenation}, \ref{LemmaConverg} we can apply Theorems 9, 10
from \cite{KKV2014} and obtain that%
\begin{equation}
\mathcal{A=}\{\gamma\left(  0\right)  :\gamma\in\mathbb{K}\}, \label{Str}%
\end{equation}
where is $\mathbb{K}$ the set of all bounded complete trajectories.

Finally, in a similar way as in \cite{KKVZ2014} let us prove that the global
attractor is stable, which means that for any $\varepsilon>0$ there is
$\delta>0$ such that%
\begin{equation}
G(t,O_{\delta}(\mathcal{A)})\subset O_{\varepsilon}(\mathcal{A})\text{ for all
}t\geq0\text{,} \label{Stable}%
\end{equation}
where $O_{\eta}(\mathcal{A})=\{z\in H:dist(z,\mathcal{A)<\eta}\}.$

\begin{lemma}
Assume that $\beta\geq3$. The global attractor $\mathcal{A}$ given in Theorem
\ref{AttrExist} is stable.
\end{lemma}

\begin{proof}
By contradiction if (\ref{Stable}) does not hold, then there exist
$\varepsilon>0$ and sequences $\delta_{k}\rightarrow0$, $x_{k}\in
O_{\delta_{k}}(\mathcal{A)}$, $t_{k}\geq0,\ y_{k}\in G(t_{k},x_{k})$ such
that
\begin{equation}
dist\left(  y_{k},\mathcal{A}\right)  \geq\varepsilon. \label{NotStable}%
\end{equation}

We consider two cases:\ 1) $t_{k}\rightarrow+\infty$ for some subsequence; 2)
$t_{k}\leq C.$

In the first situation, as the sequence $\{x_{k}\}$ belongs to a bounded set,
by the definition of global attractor we get that $dist\left(  y_{k}%
,\mathcal{A}\right)  \rightarrow0$, which contradicts (\ref{NotStable}).

In the second one, up to a subsequence $t_{k}\rightarrow t_{0},\ x_{k}%
\rightarrow x_{0}\in\mathcal{A}$, so by Lemma \ref{LemmaConverg} and the
invariance of $\mathcal{A}$ we obtain that
\[
y_{k}\rightarrow y\in G(t_{0},x_{0})\subset G(t_{0},\mathcal{A})\subset
\mathcal{A},
\]
which is again a contradiction.
\end{proof}

\subsection{Case 2:$\ \beta>3$ or $\beta=3,\ 4\alpha\mu\geq1$}

In view of Theorem \ref{UniquenessWeak} we can define the semigroup of
operators $S:\mathbb{R}^{+}\times H\rightarrow H$ by%
\[
S(t,u_{0})=u\left(  t\right)  ,
\]
where $u\left(  \text{\textperiodcentered}\right)  $ is the unique solution to
problem (\ref{NS}) with initial condition $u_{0}$. It is straightforward to
see that $S$ satisfies the semigroup properties:\ $S\left(  0,u_{0}\right)
=u_{0}$, for any $u_{0}\in H$, and $S\left(  t+s,u_{0}\right)  =S\left(
t,S\left(  s,u_{0}\right)  \right)  $, for any $u_{0}\in H$, $t,s\geq0$. Also,
making use again of Theorem \ref{UniquenessWeak} we obtain that $S(t,u_{0})$
is continuous with respect to the initial condition $u_{0}$ for fixed $t\geq0$.

We recall that the set $\mathcal{A}$ is said to be a global attractor for $S$
if it is invariant, i.e. $S\left(  t,\mathcal{A}\right)  =\mathcal{A}$, for
all $t\geq0$, and it attracts every bounded subset $B$ of the phase space $H$,
which means that%
\[
dist_{H}\left(  S\left(  t,B\right)  ,\mathcal{A}\right)  \rightarrow0\text{
as }t\rightarrow+\infty\text{,}%
\]
where $dist_{X}\left(  C,A\right)  =\sup_{x\in C}\inf_{y\in A}\left\Vert
x-y\right\Vert _{X}$ is the Hausdorff semidistance between subsets of the
Banach space $X$.

Usually in the literature a global attractor is supposed to be compact as
well. However, we prefer to use this more general definition and add
compactness as an additional property, as generally speaking a global
attractor does not have to be bounded (see \cite{Valero1999} for a non-trivial
example of an unbounded non-locally compact attractor).

The existence of the global compact attractor follows directly from Theorem
\ref{AttrExist} as a particular case. Nevertheless, we will explain this
result using the theory of attractors for semigroups as well.

For a semigroup $S$ the concepts of absorbing set and asymptotically
compactness are given in exactly the same way as for semiflows. It follows
from (\ref{EstWeak1}) that the ball
\[
B_{0}=\left\{  u\in H:\left\vert u\right\vert ^{2}\leq1+\frac{\left\lceil
f\right\rceil ^{2}}{\mu^{2}\lambda_{1}^{2}}\right\}
\]
is absorbing for the semigroup $S.$ Also, Lemma \ref{LemmaConverg} implies,
arguing as in the proof of Theorem \ref{AttrExist}, that the semigroup $S$ is
asymptotically compact.

The existence of a bounded absorbing set and the asymptotic compactness ensure
the existence of the global compact attractor \cite{LadBook}. Also, as the
space $H$ is connected, the attractor is connected \cite[p.4]{Gobbino}. Hence,
we have obtained the following result.

\begin{theorem}
\label{ExistWeakAttr}If $\beta>3$ or $\beta=3,\ 4\alpha\mu\geq1$, the
semigroup $S$ possesses the global compact connected attractor $\mathcal{A}$.
\end{theorem}

\bigskip

It is possible to prove that the global attractor is more regular if $\beta
\in(3,5)$ or $\beta=3$ and $\mu a>1/4$. Indeed, let us check that
$\mathcal{A}$ is in fact bounded in the space $\left(  H^{2}\left(
\Omega\right)  \right)  ^{3}$, and then compact in $V$ and $\left(
L^{\beta+1}\left(  \Omega\right)  \right)  ^{3}.$

\begin{lemma}
\label{EstWeakSol2}Let $\beta\in(3,5)$ or $\beta=3$ and $\mu a>1/4$. Then any
weak solution of (\ref{NS}) with initial data such that $\left\vert
u_{0}\right\vert \leq R$ satisfies the estimate%
\begin{equation}
\left\vert u_{t}(\overline{t}+r)\right\vert ^{2}+\left\Vert u\left(
\overline{t}+r\right)  \right\Vert ^{2}+\left\vert u\left(  \overline
{t}+r\right)  \right\vert _{\beta+1}^{\beta+1}\leq D\left(  R,r\right)
,\label{EstNorms}%
\end{equation}
for any $r>0$ and $\overline{t}\geq0$, where $D\left(  R,r\right)  $ is such
that $D\left(  R,r\right)  \rightarrow\infty$ if $r\rightarrow0^{+}$ or
$R\rightarrow+\infty$.
\end{lemma}

\begin{proof}
Since we have uniqueness of the Cauchy problem, the following formal
calculations can be justified via Galerkin Approximations.

We prove first the result for $\overline{t}=0$. Multiplying the equation by
$Au$ we obtain that%
\[
\frac{1}{2}\frac{d}{dt}\left\Vert u\right\Vert ^{2}+\mu\left\vert
Au\right\vert ^{2}=-b(u,u,Au)-\alpha\left(  \left\vert u\right\vert ^{\beta
-1}u,Au\right)  +\left(  f,Au\right)  .
\]
By%
\[
\left\vert b(u,u,Au)\right\vert \leq C_{1}\left\Vert u\right\Vert ^{\frac
{3}{2}}\left\vert Au\right\vert ^{\frac{3}{2}}\leq\frac{\mu}{8}\left\vert
Au\right\vert ^{2}+C_{2}\left\Vert u\right\Vert ^{6},
\]%
\[
\left\vert \left(  f,Au\right)  \right\vert \leq\frac{\mu}{4}\left\vert
Au\right\vert ^{2}+\frac{1}{\mu}\left\vert f\right\vert ^{2},
\]%
\[
\alpha\left\vert \left(  \left\vert u\right\vert ^{\beta-1}u,Au\right)
\right\vert \leq\frac{\mu}{8}\left\vert Au\right\vert ^{2}+C_{3}\left\vert
u\right\vert _{2\beta}^{2\beta},
\]
where we have used inequality (9,27) in \cite{Robinson}, we get%
\[
\frac{1}{2}\frac{d}{dt}\left\Vert u\right\Vert ^{2}+\frac{\mu}{2}\left\vert
Au\right\vert ^{2}\leq\frac{1}{\mu}\left\vert f\right\vert ^{2}+C_{2}%
\left\Vert u\right\Vert ^{6}+C_{3}\left\vert u\right\vert _{2\beta}^{2\beta}.
\]
Also, $\left\vert u\right\vert _{\infty}^{2}\leq C_{4}\left\Vert u\right\Vert
\left\vert Au\right\vert $ \cite{Temam83} gives%
\begin{align*}
\left\vert u\right\vert _{2\beta}^{2\beta}  &  =\int_{\Omega}\left\vert
u\right\vert ^{2\beta-6}\left\vert u\right\vert ^{6}dx\leq\left\vert
u\right\vert _{\infty}^{2\beta-6}\left\vert u\right\vert _{6}^{6}\\
&  \leq C_{4}\left\Vert u\right\Vert ^{\beta+3}\left\vert Au\right\vert
^{\beta-3}\leq\frac{\mu}{4C_{3}}\left\vert Au\right\vert ^{2}+C_{5}\left\Vert
u\right\Vert ^{\frac{2(\beta+3)}{5-\beta}}.
\end{align*}
Thus,%
\begin{equation}
\frac{d}{dt}\left\Vert u\right\Vert ^{2}+\frac{\mu}{2}\left\vert Au\right\vert
^{2}\leq C_{6}(1+\left\Vert u\right\Vert ^{2q}), \label{IneqQ}%
\end{equation}
where $q=(\beta+3)/(5-\beta)$. The function $y(t)=1+\left\Vert u(t)\right\Vert
^{2}$ satisfies%
\[
y^{\prime}\leq C_{6}(1+(y(t)-1)^{q})\leq C_{6}y^{q},
\]
where we have used that $1+(y-1)^{q}\leq y^{q}$ for any $y\geq1.$ Hence,%
\[
\left(  y\left(  t\right)  \right)  ^{1-q}\geq((y(t_{0}))^{1-q}+\left(
1-q\right)  C_{6}(t-t_{0}),
\]%
\[
y(t)\leq\frac{y(t_{0})}{\left(  1+(y(t_{0}))^{q-1}\left(  1-q\right)
C_{6}(t-t_{0})\right)  ^{1/(q-1)}},
\]
which is finite if%
\[
t<t_{0}+\frac{1}{(y(t_{0}))^{q-1}(q-1)C_{6}}.
\]
Take $T^{\prime}=t_{0}+1/(2(y(t_{0}))^{q-1}(q-1)C_{6})$. Then%
\begin{equation}
\left\Vert u(t)\right\Vert ^{2}\leq2^{1/(q-1)}\left(  1+\left\Vert
u(t_{0})\right\Vert ^{2}\right)  \text{ for }t\in\lbrack t_{0},T^{\prime}].
\label{IneqnormH1}%
\end{equation}
By Lemma \ref{EstWeakSol} we have%
\[
\mu\int_{0}^{r}\left\Vert u(s)\right\Vert ^{2}ds\leq\left\vert u_{0}%
\right\vert ^{2}+\frac{1+r\mu\lambda_{1}}{\mu^{2}\lambda_{1}^{2}}\left\vert
f\right\vert ^{2},
\]
for $r>0$, which implies the existence of $t_{0}\in\left(  0,r\right)  $ such
that%
\[
\left\Vert u(t_{0})\right\Vert ^{2}\leq\frac{1}{\mu r}\left(  \left\vert
u_{0}\right\vert ^{2}+\frac{1+r\mu\lambda_{1}}{\mu^{2}\lambda_{1}^{2}%
}\left\vert f\right\vert ^{2}\right)  \leq D_{1}(R,r).
\]
Hence, (\ref{IneqnormH1}) implies%
\[
\left\Vert u(t)\right\Vert ^{2}\leq2^{1/(q-1)}\left(  1+\frac{1}{\mu r}\left(
\left\vert u_{0}\right\vert ^{2}+\frac{1+r\mu\lambda_{1}}{\mu^{2}\lambda
_{1}^{2}}\left\vert f\right\vert ^{2}\right)  \right)  \leq D_{2}(R,r)\text{
}\forall\ t_{0}\leq t\leq T^{\prime}\text{.}%
\]
Let $\gamma=\min\{T^{\prime},r\}$. From (\ref{IneqQ}) we have%
\begin{equation}
\sup_{t_{0}\leq t\leq\gamma}\left\Vert u(t)\right\Vert ^{2}+\frac{\mu}{2}%
\int_{t_{0}}^{\gamma}\left\vert Au\right\vert ^{2}dt\leq D_{1}(R,r)+C_{6}%
(1+\left(  D_{2}(R,r)\right)  ^{2q})r=D_{3}(R,r). \label{IneqAu}%
\end{equation}

Multiplying the equation by $u_{t}$ and using again inequality (9.27) in
\cite{Robinson} we obtain%
\begin{align}
&  \left\vert u_{t}\right\vert ^{2}+\frac{\mu}{2}\frac{d}{dt}\left\Vert
u\right\Vert ^{2}+\frac{\alpha}{\beta+1}\frac{d}{dt}\left\vert u\right\vert
_{\beta+1}^{\beta+1}\nonumber\\
&  \leq-b(u,u,u_{t})+(f,u_{t})\leq C_{7}\left\Vert u\right\Vert ^{\frac{3}{2}%
}\left\vert Au\right\vert ^{\frac{1}{2}}\left\vert u_{t}\right\vert
+\left\vert f\right\vert ^{2}+\frac{1}{4}\left\vert u_{t}\right\vert
^{2}\nonumber\\
&  \leq\frac{1}{2}\left\vert u_{t}\right\vert ^{2}+C_{8}(\left\vert
f\right\vert ^{2}+\left\vert Au\right\vert ^{2}+\left\Vert u\right\Vert ^{6}).
\label{Inequt}%
\end{align}
Integrating over $\left(  t_{0},\gamma\right)  $ and using (\ref{IneqAu}) and
the embedding $V\subset\left(  L^{\beta+1}\left(  \Omega\right)  \right)
^{3}$ we have%
\begin{align*}
\int_{t_{0}}^{\gamma}\left\vert u_{t}\right\vert ^{2}dt  &  \leq\mu\left\Vert
u(t_{0})\right\Vert ^{2}+\frac{2\alpha}{\beta+1}\left\vert u(t_{0})\right\vert
_{\beta+1}^{\beta+1}+2C_{8}\left(  \left\vert f\right\vert ^{2}r+\frac{2}{\mu
}D_{3}(R,r)+\left(  D_{3}(R,r)\right)  ^{3}r\right) \\
&  \leq\mu D_{1}(R,r)+C_{9}\left(  D_{1}(R,r)\right)  ^{\frac{\beta+1}{2}%
}+2C_{8}\left(  \left\vert f\right\vert ^{2}r+\frac{2}{\mu}D_{3}(R,r)+\left(
D_{3}(R,r)\right)  ^{3}r\right)  =D_{4}(R,r).
\end{align*}
Hence, there is $t_{1}\in\left(  t_{0},\gamma\right)  $ such that%
\begin{equation}
\left\vert u_{t}(t_{1})\right\vert ^{2}\leq\frac{D_{4}(R,r)}{\gamma-t_{0}}%
\leq\frac{D_{4}(R,r)}{T^{\prime}-t_{0}}\leq D_{4}(R,r)2((1+D_{1}%
(R,r))^{q-1}(q-1)C_{6})=D_{5}(R,r). \label{Estut}%
\end{equation}

Further, we differentiate the equation with respect to $t$ and multiply by
$u_{t}:$%
\begin{align*}
&  \frac{1}{2}\frac{d}{dt}\left\vert u_{t}\right\vert ^{2}+\mu\left\Vert
u_{t}\right\Vert ^{2}+\alpha(\left(  \left\vert u\right\vert ^{\beta
-1}u\right)  _{t},u_{t})\\
&  =-b(u_{t},u,u_{t})-b(u,u_{t},u_{t})=b(u_{t},u_{t},u)\leq\varepsilon_{1}%
\mu\left\Vert u_{t}\right\Vert ^{2}+\frac{1}{4\varepsilon_{1}\mu}\left\vert
\left\vert u\right\vert \left\vert u_{t}\right\vert \right\vert ^{2},
\end{align*}
for $\varepsilon_{1}>0.$ Using%
\[
(\left(  \left\vert u\right\vert ^{\beta-1}u\right)  _{t},u_{t})=(\left\vert
u\right\vert ^{\beta-1}u_{t},u_{t})+\left(  \beta-1\right)  \int_{\Omega
}\left\vert u\right\vert ^{\beta-1}\left\vert u_{t}\right\vert ^{2}%
dx\geq(\left\vert u\right\vert ^{\beta-1}u_{t},u_{t})
\]
we obtain%
\begin{equation}
\frac{1}{2}\frac{d}{dt}\left\vert u_{t}\right\vert ^{2}+\mu(1-\varepsilon
_{1})\left\Vert u_{t}\right\Vert ^{2}+\alpha\left\vert \left\vert u\right\vert
^{\frac{\beta-1}{2}}\left\vert u_{t}\right\vert \right\vert ^{2}\leq\frac
{1}{4\varepsilon_{1}\mu}\left\vert \left\vert u\right\vert \left\vert
u_{t}\right\vert \right\vert ^{2}. \label{Inequt2}%
\end{equation}
From (\ref{Inequt}) and%
\[
-b(u,u,u_{t})=b(u,u_{t},u)\leq\varepsilon_{2}\mu\left\Vert u_{t}\right\Vert
^{2}+\frac{1}{4\mu\varepsilon_{2}}\left\vert u\right\vert _{4}^{4},
\]
for $\varepsilon_{2}>0$, we have%
\begin{equation}
\frac{1}{2}\left\vert u_{t}\right\vert ^{2}+\frac{\mu}{2}\frac{d}%
{dt}\left\Vert u\right\Vert ^{2}+\frac{\alpha}{\beta+1}\frac{d}{dt}\left\vert
u\right\vert _{\beta+1}^{\beta+1}\leq\varepsilon_{2}\mu\left\Vert
u_{t}\right\Vert ^{2}+\frac{1}{4\mu\varepsilon_{2}}\left\vert u\right\vert
_{4}^{4}+\frac{1}{2}\left\vert f\right\vert ^{2}. \label{Inequt3}%
\end{equation}
Summing (\ref{Inequt2}) and (\ref{Inequt3}) we get%
\begin{align*}
&  \frac{d}{dt}\left(  \frac{1}{2}\left\vert u_{t}\right\vert ^{2}+\frac{\mu
}{2}\left\Vert u\right\Vert ^{2}+\frac{\alpha}{\beta+1}\left\vert u\right\vert
_{\beta+1}^{\beta+1}\right)  +\frac{1}{2}\left\vert u_{t}\right\vert ^{2}%
+\mu(1-\varepsilon_{1}-\varepsilon_{2})\left\Vert u_{t}\right\Vert ^{2}%
+\alpha\left\vert \left\vert u\right\vert ^{\frac{\beta-1}{2}}\left\vert
u_{t}\right\vert \right\vert ^{2}\\
&  \leq\frac{1}{4\varepsilon_{1}\mu}\left\vert \left\vert u\right\vert
\left\vert u_{t}\right\vert \right\vert ^{2}+\frac{1}{4\mu\varepsilon_{2}%
}\left\vert u\right\vert _{4}^{4}+\frac{1}{2}\left\vert f\right\vert ^{2}%
\leq\frac{1}{4\varepsilon_{1}\mu}\left\vert \left\vert u\right\vert \left\vert
u_{t}\right\vert \right\vert ^{2}+K_{1}(\varepsilon_{2})\left(  1+\left\vert
u\right\vert _{\beta+1}^{\beta+1}\right)  .
\end{align*}
If $\beta=3$, the condition $\alpha\mu4>1$ implies that
\[
\alpha\left\vert \left\vert u\right\vert ^{\frac{\beta-1}{2}}\left\vert
u_{t}\right\vert \right\vert ^{2}-\frac{1}{4\varepsilon_{1}\mu}\left\vert
\left\vert u\right\vert \left\vert u_{t}\right\vert \right\vert ^{2}=\left(
\alpha-\frac{1}{4\varepsilon_{1}\mu}\right)  \left\vert \left\vert
u\right\vert \left\vert u_{t}\right\vert \right\vert ^{2}\geq0
\]
for $\varepsilon_{1}$ close enough to $1$. If $\beta>3$, then by Young's
inequality there is $K_{2}(\varepsilon_{1})>0$ such that%
\[
\alpha\left\vert \left\vert u\right\vert ^{\frac{\beta-1}{2}}\left\vert
u_{t}\right\vert \right\vert ^{2}-\frac{1}{4\varepsilon_{1}\mu}\left\vert
\left\vert u\right\vert \left\vert u_{t}\right\vert \right\vert ^{2}\geq
\frac{\alpha}{2}\left\vert \left\vert u\right\vert ^{\frac{\beta-1}{2}%
}\left\vert u_{t}\right\vert \right\vert ^{2}-K_{2}(\varepsilon_{1})\left\vert
u_{t}\right\vert ^{2}.
\]
Hence,%
\[
\frac{d}{dt}\left(  \frac{1}{2}\left\vert u_{t}\right\vert ^{2}+\frac{\mu}%
{2}\left\Vert u\right\Vert ^{2}+\frac{\alpha}{\beta+1}\left\vert u\right\vert
_{\beta+1}^{\beta+1}\right)  \leq K_{1}(\varepsilon_{2})\left(  1+\left\vert
u\right\vert _{\beta+1}^{\beta+1}\right)  +K_{2}(\varepsilon_{1})\left\vert
u_{t}\right\vert ^{2},\ \forall\ t\geq t_{1},
\]
if we choose $\varepsilon_{1}$ close enough to $1$ and $\varepsilon_{2}$ small
enough. Let $y(t)=\frac{1}{2}\left\vert u_{t}\right\vert ^{2}+\frac{\mu}%
{2}\left\Vert u\right\Vert ^{2}+\frac{\alpha}{\beta+1}\left\vert u\right\vert
_{\beta+1}^{\beta+1}$. Then using Gronwall lemma, the embedding $V\subset
\left(  L^{\beta+1}(\Omega)\right)  ^{3}$ and (\ref{IneqAu}), (\ref{Estut}) we
obtain for some constants $K=K(\varepsilon_{1},\varepsilon_{2}),\ D_{6}(R,r)$
that%
\[
y(t)\leq\left(  y(t_{1})+1\right)  e^{K(t-t_{1})}\leq D_{6}(R,r)e^{K(t-t_{1}%
)}\ \forall t\geq t_{1}.
\]
Thus, in particular,%
\[
y(r)\leq D_{6}(R,r)e^{K(r-t_{1})}=D_{7}(R,r),
\]
which proves (\ref{EstNorms}) for $\overline{t}=0.$

For an arbitary $\overline{t}\geq0$ by Lemma \ref{EstWeakSol} we make use of
the estimate%
\[
\left\vert u(t)\right\vert ^{2}\leq\left\vert u(0)\right\vert ^{2}%
+\frac{\left\vert f\right\vert ^{2}}{\mu^{2}\lambda_{1}^{2}}\leq R^{2}%
+\frac{\left\vert f\right\vert ^{2}}{\mu^{2}\lambda_{1}^{2}}=R_{1}%
^{2},\ \forall t\geq0.
\]
Defining $v(t)=u(t+\overline{t})$, then
\[
\frac{1}{2}\left\vert v_{t}\left(  r\right)  \right\vert ^{2}+\frac{\mu}%
{2}\left\Vert v(r)\right\Vert ^{2}+\frac{\alpha}{\beta+1}\left\vert v\left(
r\right)  \right\vert _{\beta+1}^{\beta+1}\leq D_{7}\left(  R_{1},r\right)
=D_{8}(R,r),
\]
which gives (\ref{EstNorms}).
\end{proof}

\bigskip

As a consequence of this lemma and the compact embedding $V\subset H$ we
obtain the following result.

\begin{corollary}
\label{CompactOperator}For any $r>0$ the map $u_{0}\mapsto S\left(
r,u_{0}\right)  $ maps bounded subsets of $H$ onto bounded subsets of $V\cap
L^{\beta+1}\left(  \Omega\right)  $. Hence, $S\left(  r\right)  $ is a compact
operator, i.e., it maps bounded subsets of $H$ onto relatively compact ones.
\end{corollary}

\begin{lemma}
\label{EstAu}Let $\beta\in(3,5)$ or $\beta=3$ and $\mu a>1/4$. Then any weak
solution of (\ref{NS}) with initial data such that $\left\vert u_{0}%
\right\vert \leq R$ satisfies the estimate%
\[
\left\vert Au\left(  r\right)  \right\vert \leq K\left(  R,r\right)  ,
\]
for any $r>0$, where $K\left(  R,r\right)  $ is such that $K\left(
R,r\right)  \rightarrow\infty$ if $r\rightarrow0^{+}$ or $R\rightarrow+\infty$.
\end{lemma}

\begin{proof}
By Proposition 9.2 in \cite{Robinson} we have%
\[
\left\vert B(u,u)\right\vert \leq d_{1}\left\Vert u\right\Vert ^{\frac{3}{2}%
}\left\vert Au\right\vert ^{\frac{1}{2}}\leq\frac{\mu}{4}\left\vert
Au\right\vert +d_{2}\left\Vert u\right\Vert ^{3}.
\]
Using the Gagliardo-Nirenberg inequality and $\beta<5$ we find that%
\begin{align*}
\alpha\left\vert \left\vert u\right\vert ^{\beta-1}u\right\vert  &
=\alpha\left\vert u\right\vert _{2\beta}^{\beta}\leq d_{3}\left\vert
Au\right\vert ^{\frac{3\left(  \beta-1\right)  }{\beta+7}}\left\vert
u\right\vert _{\beta+1}^{\frac{\beta^{2}+4\beta+3}{\beta+7}}\\
&  \leq\frac{\mu}{4}\left\vert Au\right\vert +d_{4}\left\vert u\right\vert
_{\beta+1}^{\frac{\beta^{2}+4\beta+3}{10-2\beta}}.
\end{align*}

Hence,
\[
\frac{\mu}{2}\left\vert Au\left(  r\right)  \right\vert \leq\left\vert
u_{t}(r)\right\vert +d_{2}\left\Vert u(r)\right\Vert ^{3}+d_{4}\left\vert
u(r)\right\vert _{\beta+1}^{\frac{\beta^{2}+4\beta+3}{10-2\beta}}+\left\vert
f\right\vert ,
\]
so the result follows by applying Lemma \ref{EstWeakSol2}.
\end{proof}

\bigskip

We are now in position of proving the regularity of the global attractor.

\begin{theorem}
\label{AttrReg}Let $\beta\in(3,5)$ or $\beta=3$ and $\mu a>1/4$. Then the
global attractor $\mathcal{A}$ is bounded in $\left(  H^{2}(\Omega)\right)
^{3}$, and then compact in $V$ and $\left(  L^{\beta+1}(\Omega)\right)  ^{3}$.
Moreover,
\begin{equation}
dist_{V}(S(t,B),\mathcal{A})\rightarrow0\text{ as }t\rightarrow+\infty,
\label{ConvergVLbeta}%
\end{equation}
for any $B$ bounded in $H$.
\end{theorem}

\begin{proof}
Since the global attractor is invariant, $\mathcal{A}=S(r,\mathcal{A})$ for
$r>0$, so $\mathcal{A}$ is bounded in $\left(  H^{2}(\Omega)\right)  ^{3}$ by
Lemma \ref{EstAu}. The compact embeddings $H^{2}(\Omega)\subset H^{1}(\Omega
)$, $H^{2}(\Omega)\subset L^{\beta+1}(\Omega)$ imply the compactness of the
attractor in $V$ and $\left(  L^{\beta+1}(\Omega)\right)  ^{3}$.

Let $B_{0}$ be an absorbing ball. By Lemma \ref{EstAu} the set $B_{1}%
=S(r,B_{0})$ is bounded in $\left(  H^{2}(\Omega)\right)  ^{3}$ and%
\[
S(t,B)=S(r,S(t-r,B))\subset B_{1}\text{ for }t\geq t_{0}(B)\text{.}%
\]
From here it is easy to deduce (\ref{ConvergVLbeta}).
\end{proof}

\bigskip

In the paper \cite{PardoValGim} the results of Theorem \ref{EstWeakSol2},
Lemma \ref{EstAu} and Theorem \ref{AttrReg} were stated for any $\beta\geq4$.
However, as mentioned in the introduction, the proof is not correct because in
the case of bounded domains we cannot multiply the equation by $-\Delta u$ to
obtain the regularity of weak solutions. Hence, these results remain open for
$\beta\geq5.$

\bigskip

We recall that a function $u:[0,T]\rightarrow V\cap\left(  L^{\beta+1}%
(\Omega)\right)  ^{3}$ is called a strong solution of (\ref{NS}) if $u$ is a
weak solution and
\[
u\in L^{\infty}(0,T;V)\cap L^{2}(0,T;D(A))\cap L^{\infty}(0,T;(L^{\beta
+1}(\Omega))^{3}).
\]

\begin{theorem}
\label{ExistStronSol}\cite[Theorem 1.1]{KimLi} Let $u_{0}\in V\cap\left(
L^{\beta+1}(\Omega)\right)  ^{3}$, $f\in H$ and either $\beta\in(3,5)$ or
$\beta=3$ and $\mu a>1/4.$ Then there exists a unique strong solution of
(\ref{NS}).
\end{theorem}

For $\beta\in(3,5)$ or $\beta=3$ and $\mu a>1/4$ let $S_{V}:\mathbb{R}%
^{+}\times V\rightarrow V$ be the semigroup defined by the strong solutions of
(\ref{NS}).

\begin{corollary}
Let $\beta\in(3,5)$ or $\beta=3$ and $\mu a>1/4$. Then $\mathcal{A}$ is also
the global attractor for $S_{V}.$
\end{corollary}

We have given an alternative proof of the existence of the strong attractor,
as the proof given in the papers \cite{SongHou}, \cite{SongHou2} and
\cite{SongLiangWu} for $3<\beta\leq5$ is unclear. For $\beta\geq5$ the problem
remains open. Also, in \cite{PardoValGim} a conditional result was stated
about the strong attractor for $\beta>5$. For the same reason as before, the
proof is also unclear.

\bigskip

\section{Numerical simulations}

We shall now focus on solving numerically the equation (\ref{NS}) employing
computational fluid dynamics (CFD) to visualize scenarios in which the
evolution of the fluid flow converges to a steady state. It is important to
stress that the examples reported below are only intended for showing the
asymptotic behaviour of the fluid flow numerically when taking different
values of the parameters $\alpha$ and $\beta$ in the momentum equation of
(\ref{NS}), but no conclusive results should be deduced from the numerical simulations.

\begin{figure}[h]
\centering
\includegraphics[width=0.5\linewidth]{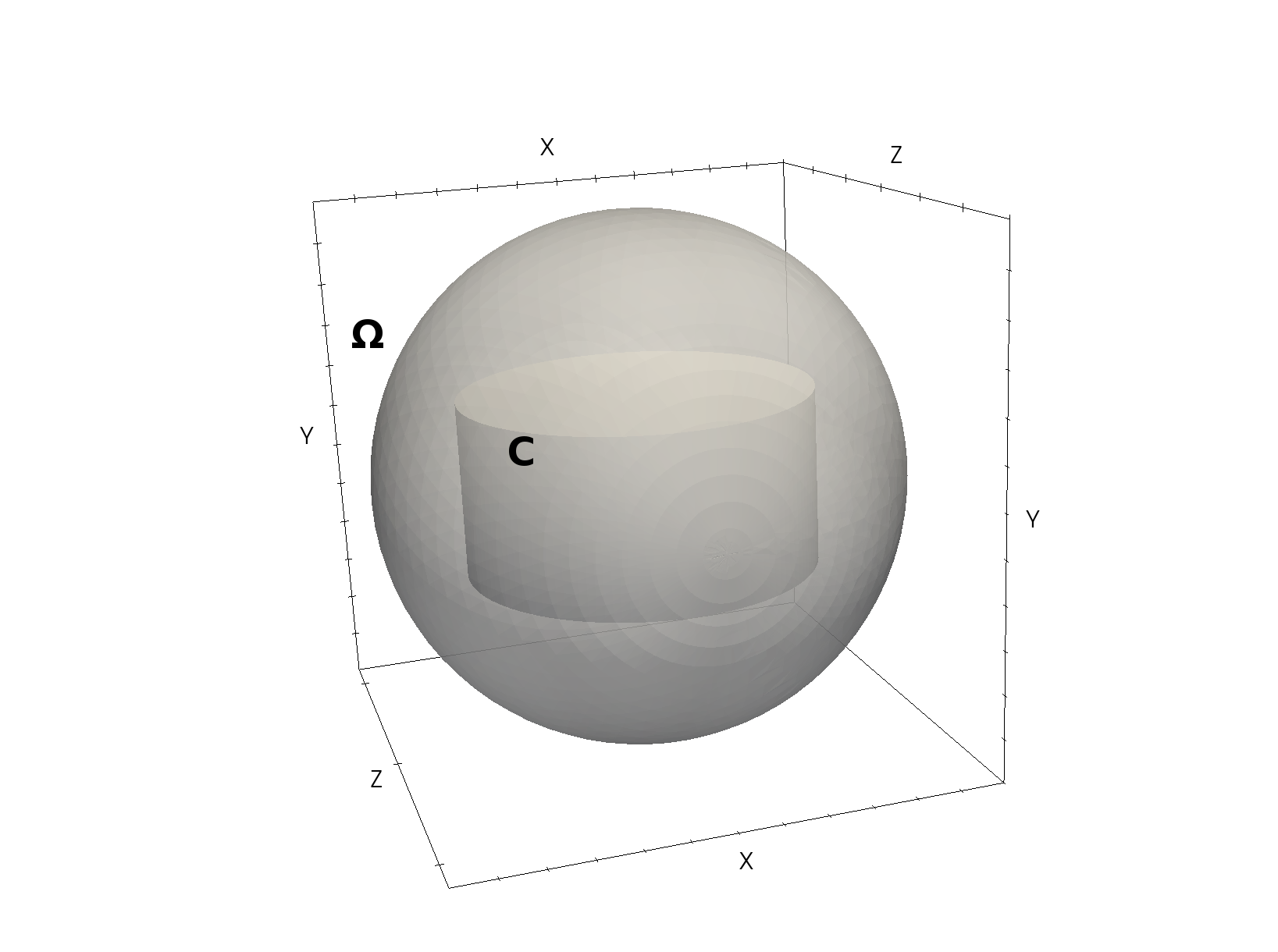} \caption{Flow domain}%
\label{flow-domain}%
\end{figure}

The geometry of the flow domain used in all our numerical experiments is a
sphere $\Omega$ of radius $6$ centered at the origin. We also take the source
term $f$ in (\ref{NS}) as
\[
f(x)=%
\begin{cases}
(0,2,0) & \text{if }x\in C\\
(0,0,0) & \text{if }x\in\Omega\setminus C
\end{cases}
\]
where $C$ is a cylinder, with both radius and height of $4$, within the flow
domain symmetrically located at the center of the sphere $\Omega$ in such a
way that the base of the cylinder is parallel to the $xz$-plane as in Figure
\ref{flow-domain}. Observe that $f(x,t)$ can be seen as a constant source
force within the cylinder $C$ propelling the fluid flow upwards.


Numerical simulations were all performed by using the CFD package
OpenFOAM\textsuperscript{\textregistered}, which is the acronym of
\textit{Open Source Field Operation and Manipulation}.
OpenFOAM\textsuperscript{\textregistered} is an open-source CFD software based
on C++ that contains a toolbox for tailored numerical solvers for a wide
variety of problems relevant to the industry and scientific community. The
solvers implemented in OpenFOAM\textsuperscript{\textregistered} uses the
Finite Volume Method (FVM) to discretize the governing equations on
unstructured meshes (see \cite{ferziger,versteeg}). The solver used to
integrate our model numerically was \emph{pimpleFoam}, which combines the two
most common algorithms for solving the Navier-Stokes equations, namely, SIMPLE
and PISO algorithms. The pimpleFoam code is inherently transient, requiring an
initial condition and boundary conditions.
OpenFOAM\textsuperscript{\textregistered} includes pre-processing and
post-processing capabilities such as snappyHexMesh and ParaFoam for meshing
and visualization, respectively.

\begin{figure}[h]
\centering
\includegraphics[width=0.45\linewidth]{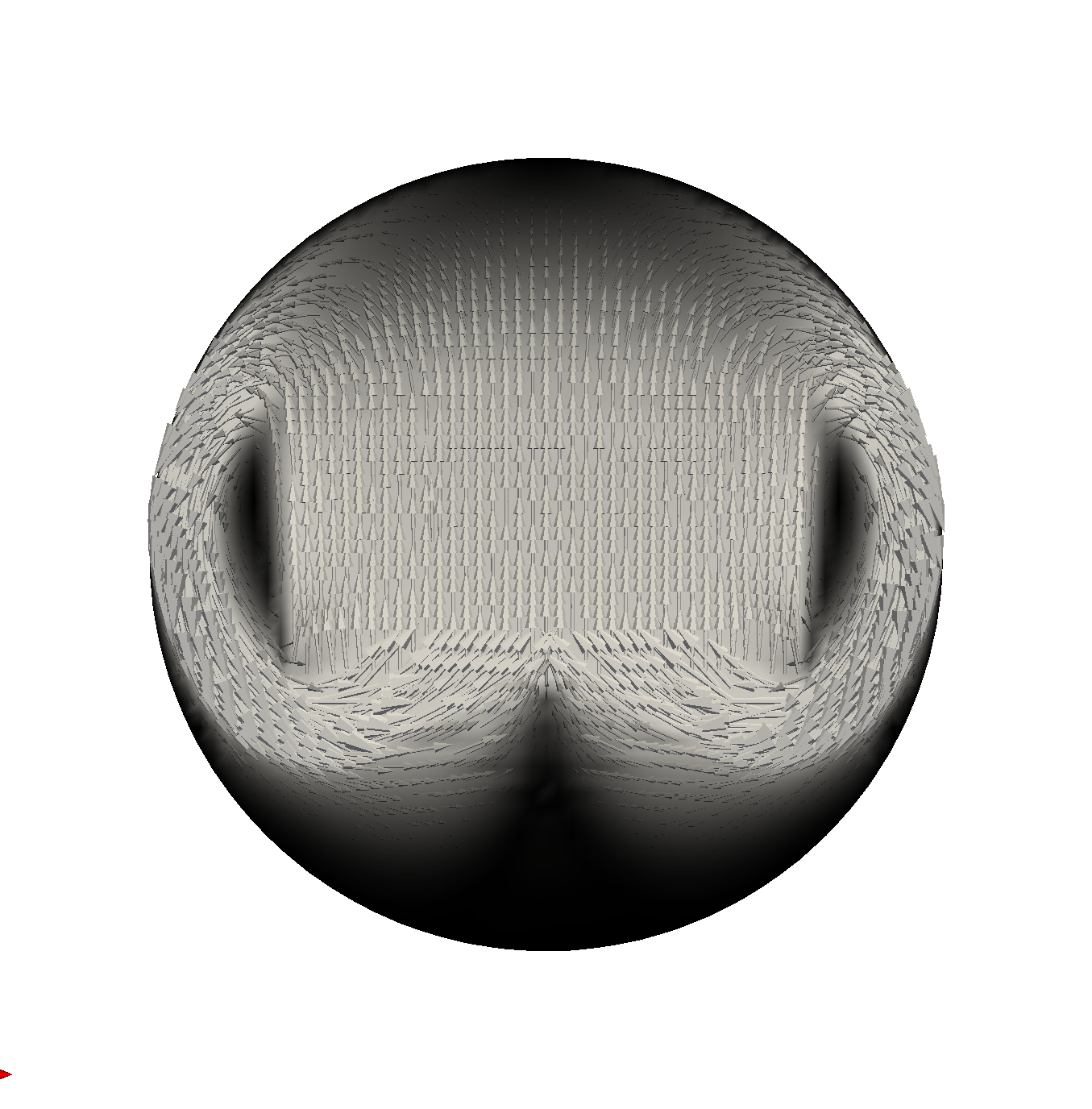}
\includegraphics[width=0.45\linewidth]{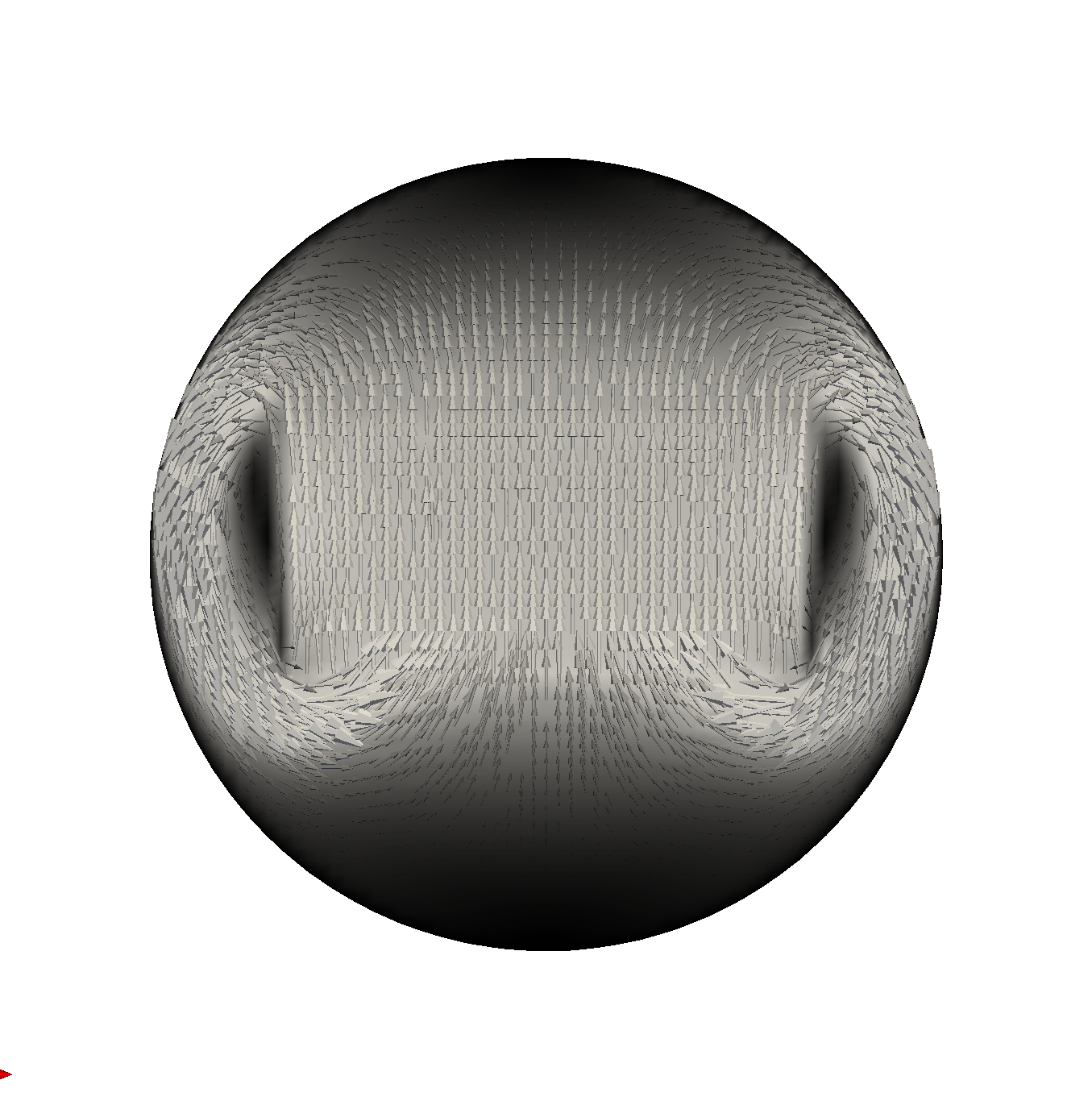} \caption{Flow velocity $u$
for the steady state in the $xy$-section at $z=0$. The darker areas mean lower
fluid flow speed. The initial velocity $u_{0}$ is identically zero. Left panel
parameters: $\alpha=0.2;\, \beta=1$. Right panel parameters: $\alpha=0.5;\,
\beta=1$. }%
\label{experimento1}%
\end{figure}

\begin{figure}[h]
\centering
\includegraphics[width=0.45\linewidth]{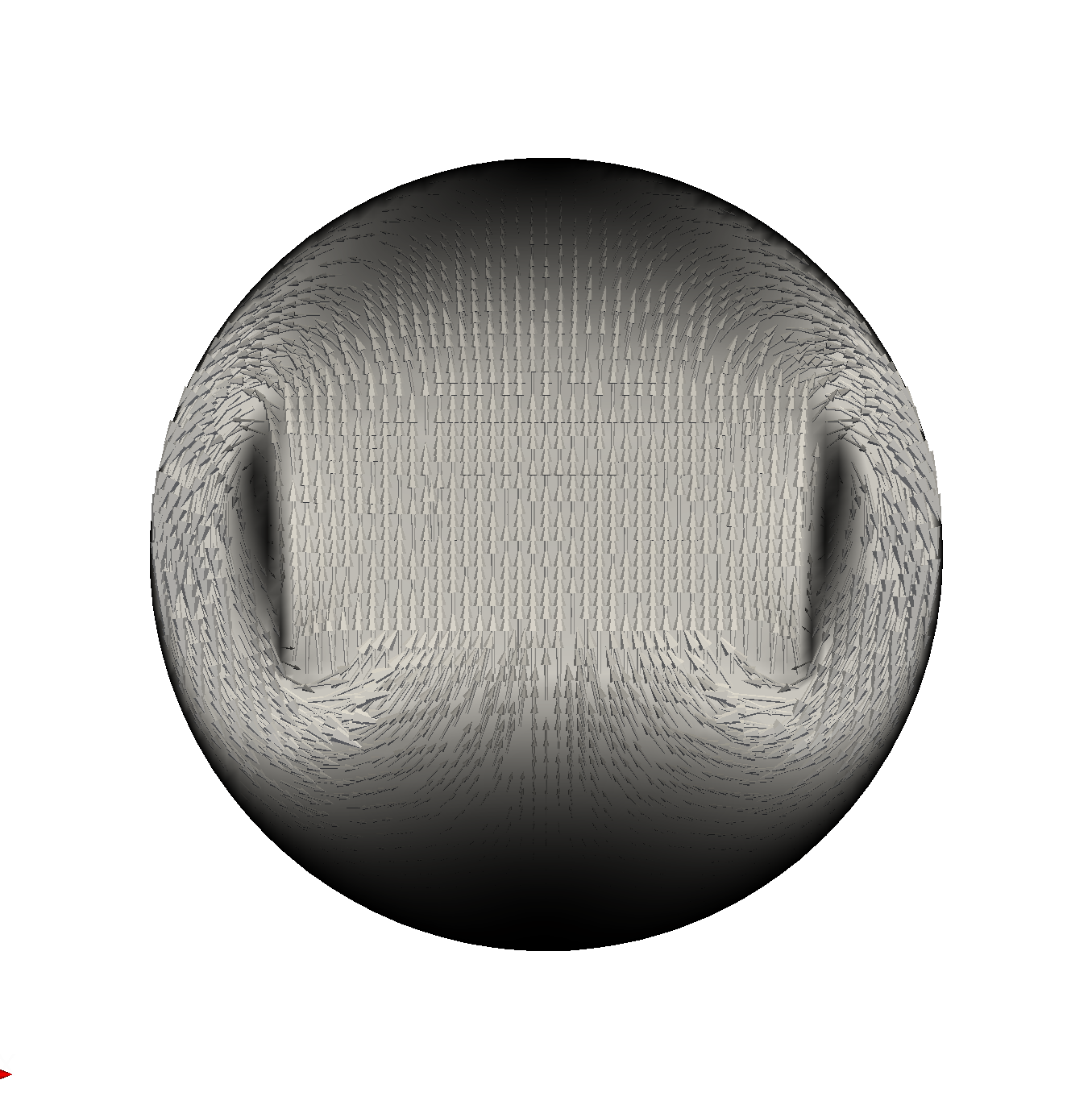}
\includegraphics[width=0.45\linewidth]{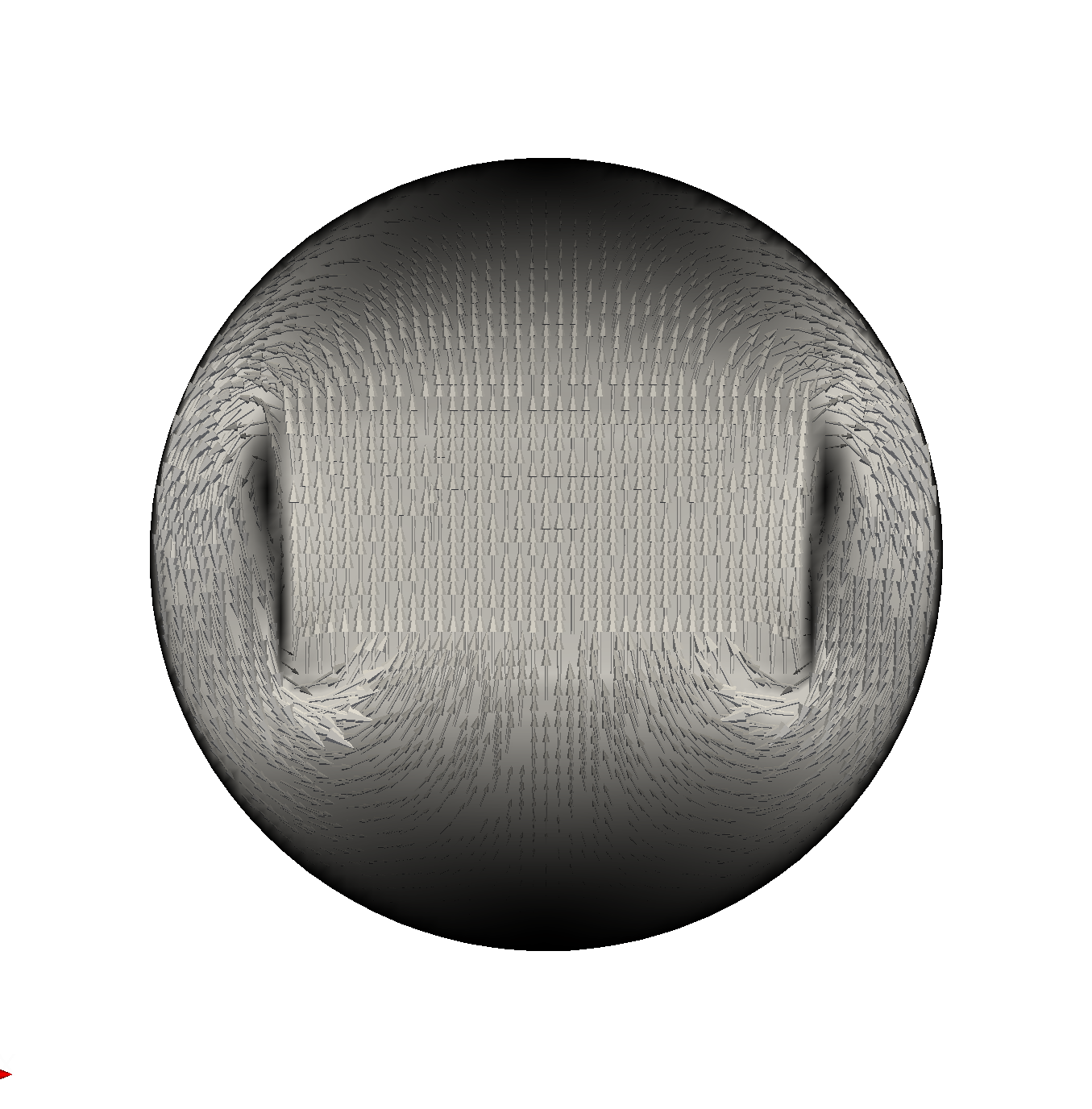} \caption{Flow velocity $u$
for the steady state in the $xy$-section at $z=0$. The darker areas mean lower
fluid flow speed. The initial velocity $u_{0}$ is identically zero. Left panel
parameters: $\alpha=0.2;\, \beta=2$. Right panel parameters: $\alpha=0.5;\,
\beta=2$.}%
\label{experimento2}%
\end{figure}

\begin{figure}[h]
\centering
\includegraphics[width=0.45\linewidth]{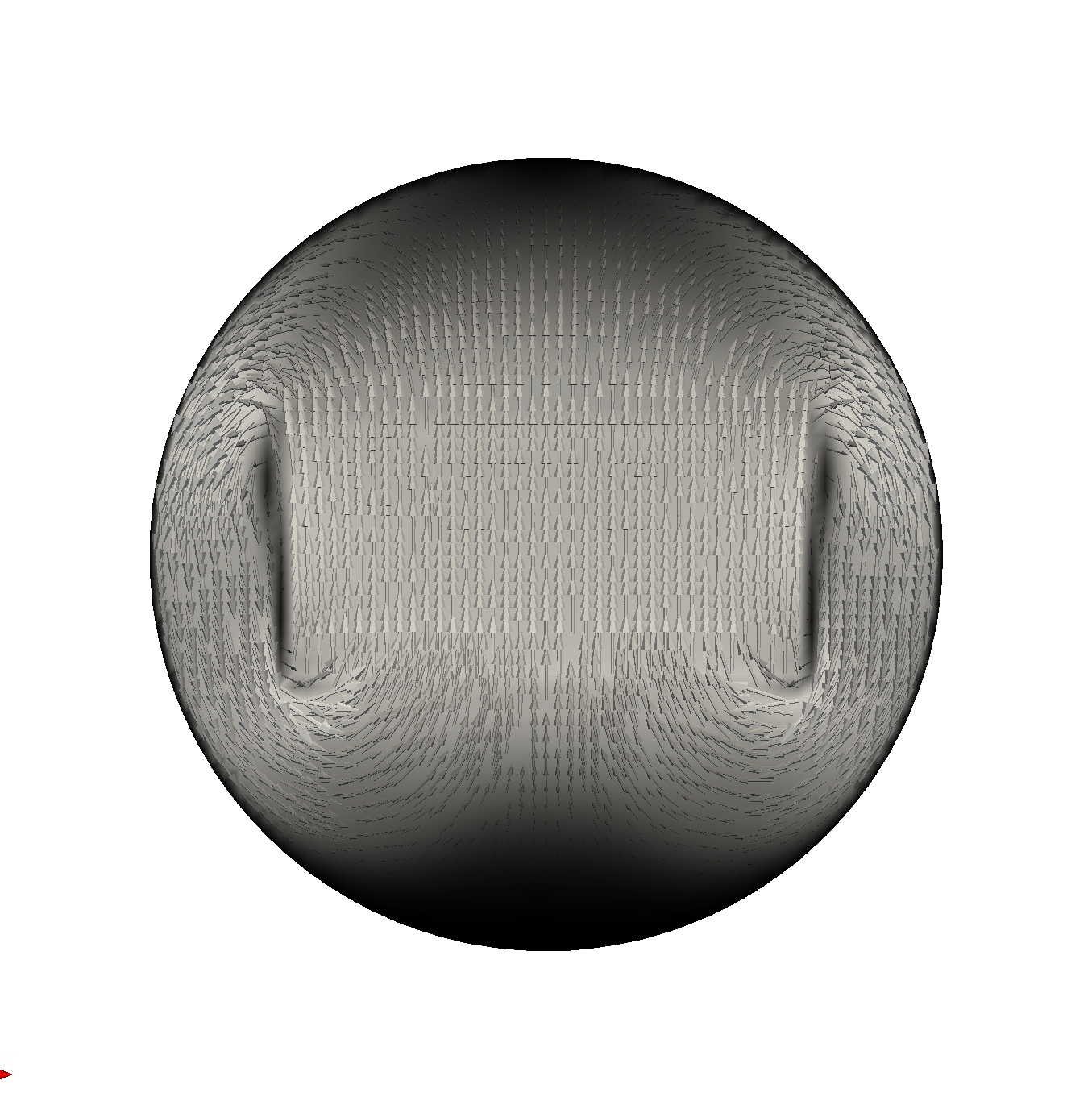}
\includegraphics[width=0.45\linewidth]{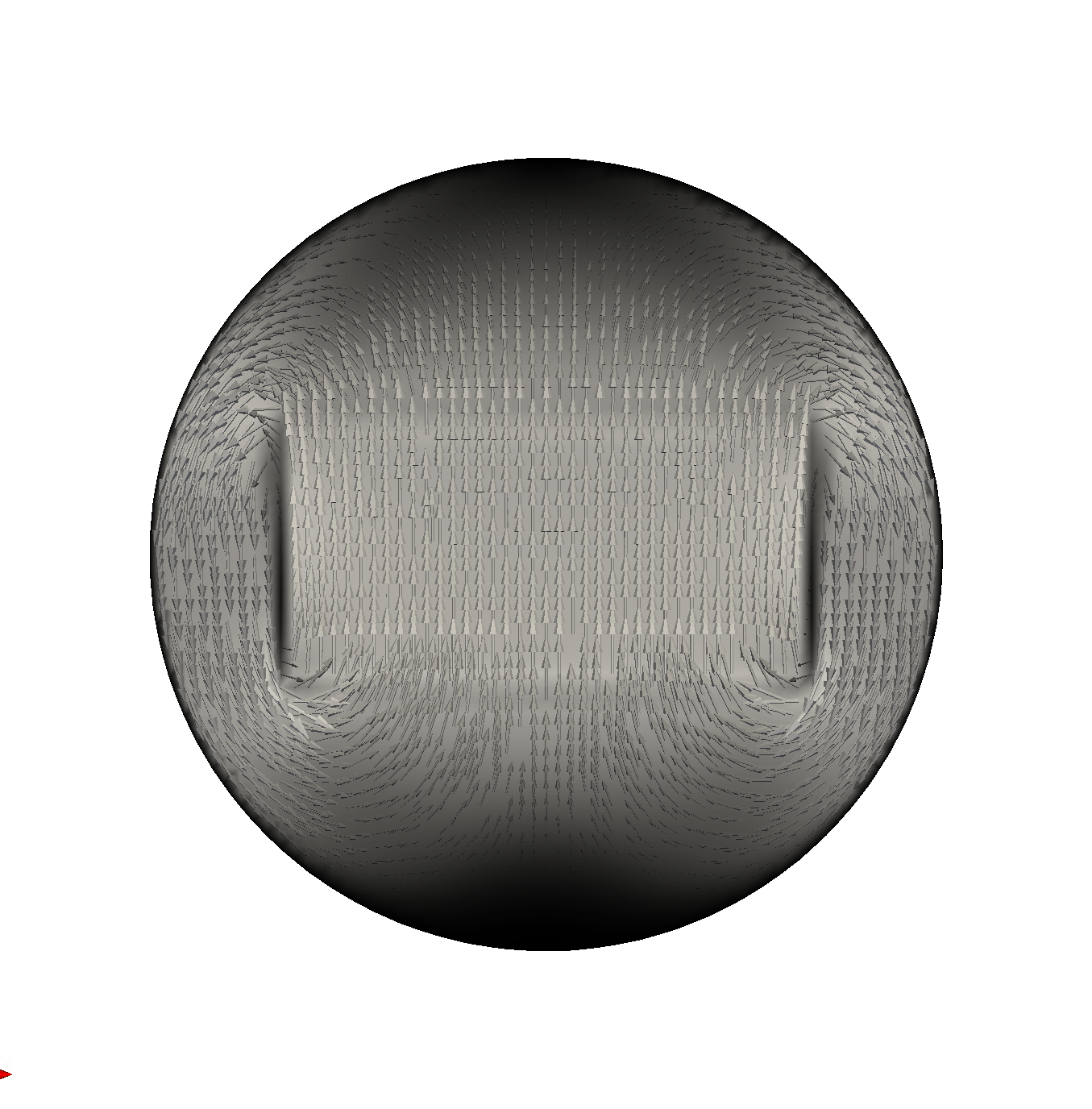}\caption{Flow velocity $u$
for the steady state in the $xy$-section at $z=0$. The darker areas mean lower
fluid flow speed. The initial velocity $u_{0}$ is identically zero. Left panel
parameters: $\alpha=0.2;\, \beta=4$. Right panel parameters: $\alpha=0.5;\,
\beta=4$.}%
\label{experimento3}%
\end{figure}

Figures \ref{experimento1}, \ref{experimento2} and \ref{experimento3} show the
steady state for the numerical solution of the equations (\ref{NS}) for
experiments with different values of the parameters $\alpha$ and $\beta$. In
those experiments we have set the initial condition $u_{0}(x)=(0,0,0)$, and
the images represent the velocity vector field in the $xy$-section at $z=0$.
The darker areas in the images are those where the magnitude for the velocity
vector $u$ is smaller, while the lighter ones represent the areas where the
velocity is higher. According to the results of the experiments, when the
magnitude of the velocity vector is greater than $1$ and the parameter $\beta$
becomes bigger, the medium provides increased resistance to movement, so the
fluid flow slows down more quickly. On the contrary, when the magnitude of the
velocity vector is less than $1$ and the parameter $\beta$ becomes smaller,
then medium provides decreased resistance to movement, so the fluid flow
spreads further through the medium. The effect of the parameter $\alpha$ does
not depend of the $u$ magnitude, acting proportionally, i.e., the larger the
value of $\alpha$, the higher resistance to motion of the fluid flow. It has
also been observed that convergence speed to the steady state is higher as
$\alpha$ and $\beta$ increase. Therefore, when $\alpha$ and $\beta$ are small,
a higher period of time to get convergence to the steady state is required. In
fact, we have also performed simulations (not shown here) for values of
$\alpha$ close to $0$ (also for $\alpha=0$) and a low value of $\beta$ (for
instance, $\beta=1$), but we did not achieve convergence to a steady state for
an approachable (from a computational point of view) time value. It is likely
that for such values of the parameters the global attractor (if it exists!) is
more complex than a fixed point.

\begin{figure}[h]
\centering
\includegraphics[width=0.45\linewidth]{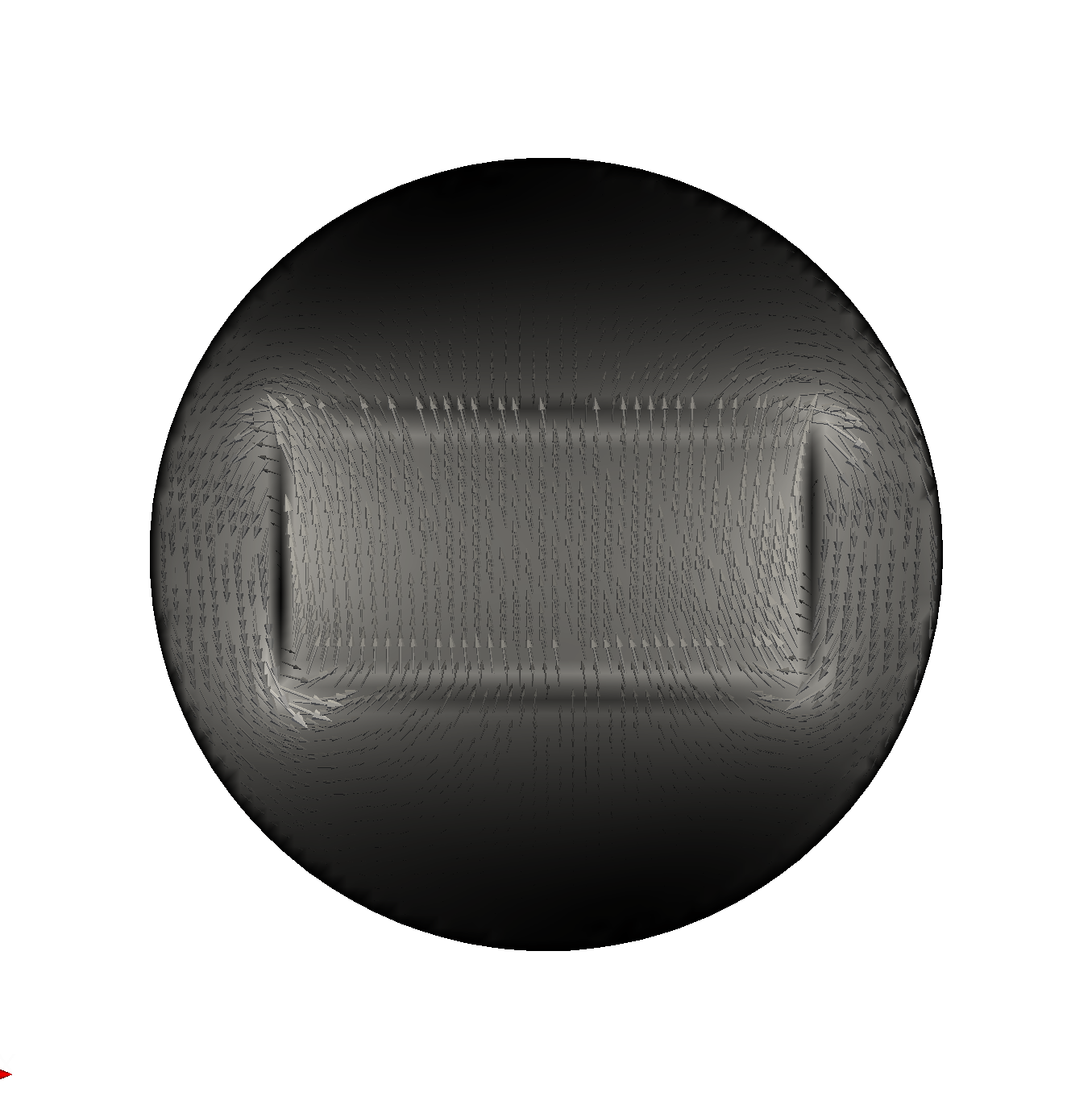}
\includegraphics[width=0.45\linewidth]{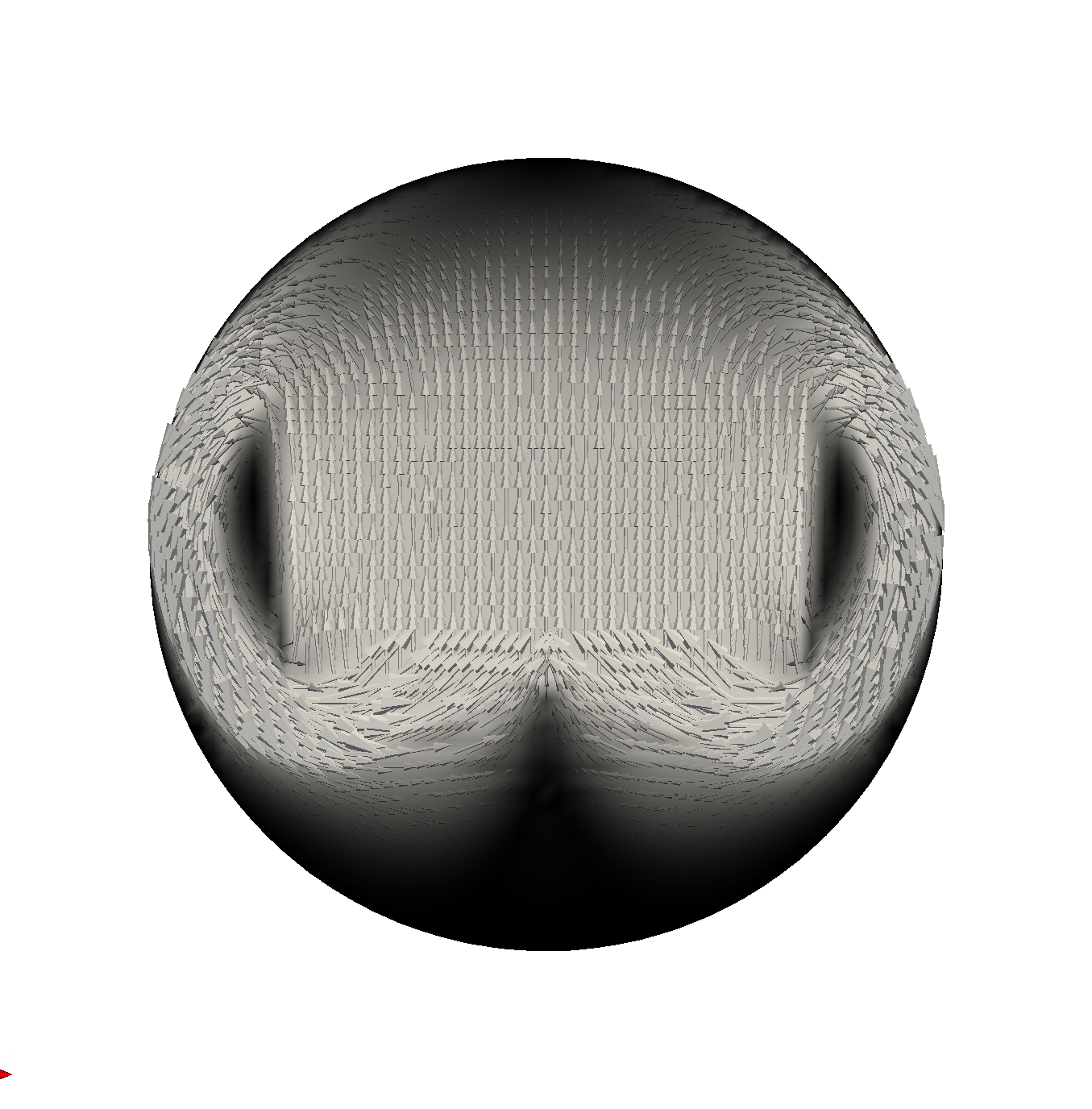} \caption{Flow velocity $u$
in the $xy$-section at $z=0$ for $\alpha=0.2$ and $\beta=1$. The darker areas
mean lower fluid flow speed. The initial velocity $u_{0}=(1,0,0)$. Left panel:
state when $t=0.1$. Right panel: steady state ($t$ large enough).}%
\label{experimento4}%
\end{figure}

Lastly, an experiment with a non-vanishing initial condition was carried out.
The performance was made by taking $u_{0}(x)=(1,0,0)$, $\alpha=0.2$ and
$\beta=1$, and the results of the model is shown in Figure \ref{experimento4}.
As one might expect, the steady state does not depend of the initial
condition, hence, it is the same as taking $u_{0}$ equal to zero (compare the
right panel of Figure \ref{experimento4} to the left panel of Figure
\ref{experimento1}).

\bigskip

\textbf{Acknowledgements}

The research has been partially supported by the Spanish Ministerio de Ciencia
e Innovaci\'{o}n (MCI), Agencia Estatal de Investigaci\'{o}n (AEI) and Fondo
Europeo de Desarrollo Regional (FEDER) under the project PID2021-122991NB-C21,
and by the Generalitat Valenciana, project PROMETEO/2021/063.

\newpage

\end{document}